\documentclass{gtart}

\def\ifplaintex{\expandafter\ifx\csname documentclass\endcsname\relax}

\def\gtp{{\mathsurround=0pt\it $\cal G\mskip-2mu$eometry \&\ 
$\cal T\!\!$opology $\cal P\!$ublications}}  

\def\recd{{\small Received:\qua\receiveddate\ifx\reviseddate\relax
\else\qquad Revised:\qua\reviseddate\fi\par}} 

\def\lognumber#1{\def\thelognumber{#1}}
\def\volumenumber#1{\def\thevolumenumber{#1}}
\def\volumeyear#1{\def\thevolumeyear{#1}}
\def\papernumber#1{\def\thepapernumber{#1}}
\def\pagenumbers#1#2{\def\startpage{#1}\def\finishpage{#2}}
\def\published#1{\def\publishdate{#1}}

\def\received#1{\def\receiveddate{#1}}

\def\accepted#1{\def\accepteddate{#1}}

\def\asciiaddress#1{\def\theasciiaddress{#1}}

\long\def\asciiabstract#1{\long\def\theasciiabstract{#1}}

\let\\\par\let\thelognumber\relax\let\thevolumenumber\relax
\let\thepapernumber\relax\let\thevolumeyear\relax\let\startpage\relax
\let\finishpage\relax\let\publishdate\relax\let\receiveddate\relax
\let\reviseddate\relax\let\accepteddate\relax\let\theasciititle\relax
\let\theasciiauthors\relax\let\theasciiaddress\relax
\let\theasciiabstract\relax

\let\theasciiemail\relax

\ifplaintex
\font\logobig=cmssbx10 scaled 3836
\font\logomed=cmssbx10 scaled 2557
\else
\font\logobig=cmssbx10 scaled 4200
\font\logomed=cmssbx10 scaled 2800
\fi

\long\def\makeagttitle{   
\count0=\startpage
\agt\hfill      
\hbox to 45truept{\vbox to 0pt{\vglue -13truept{\logomed A\kern -.37em{\logobig 
T}\kern -.38em G}\vss}\hss}
\break
{\small Volume \thevolumenumber\ (\thevolumeyear)
\startpage--\finishpage\nl
Published: \publishdate}

\vglue .25truein

{\parskip=0pt\leftskip 0pt plus
1fil\def\\{\par\smallskip}{\Large\bf\thetitle}\par\medskip} \vglue
0.05truein

{\parskip=0pt\leftskip 0pt plus 1fil\def\\{\par}{\sc\theauthors}
\par\medskip}%
 
\vglue 0.03truein 

{\small\leftskip 25truept\rightskip 25truept{\bf Abstract}\stdspace\theabstract

{\bf AMS Classification}\stdspace\theprimaryclass
\ifx\thesecondaryclass\relax\else; \thesecondaryclass\fi\par
{\bf Keywords}\stdspace \thekeywords\par}\vglue 7truept

}   

\ifplaintex
\font\phead=cmsl9 scaled 950
\font\pnum=cmbx10 scaled 913
\font\pfoot=cmsl9 scaled 950
\headline{\vbox to 0pt{\vskip -4.5mm\line{\small\phead\ifnum
\count0=\startpage ISSN 1472-2739 (on-line) 1472-2747 (printed)
\hfill {\pnum\folio}\else\ifodd\count0\def\\{ }%
\ifx\theshorttitle\relax\thetitle\else\theshorttitle\fi\hfill{\pnum\folio}
\else\def\\{ and }{\pnum\folio}\hfill\ifx\theshortauthors\relax\theauthors
\else\theshortauthors\fi\fi\fi}\vss}}
\footline{\vbox to 0pt{\vglue 0mm\line{\small\pfoot\ifnum\count0=\startpage
\copyright\ \gtp\hfill\else
\agt, Volume \thevolumenumber\ (\thevolumeyear)\hfill\fi}\vss}}
\else
\font=cmsl9 scaled 1050
\font\lnum=cmbx10 
\font=cmsl9 scaled 1050
\makeatletter
\def\@oddhead{{\small\lhead\ifnum\count0=\startpage ISSN 1472-2739 
(on-line) 1472-2747 (printed)\hfill {\lnum\number\count0}\else\ifodd\count0
\def\\{ }\ifx\theshorttitle\relax \thetitle \else\theshorttitle\fi\hfill
{\lnum\number\count0}\else\def\\{ and }{\lnum\number\count0}
\hfill\ifx\theshortauthors\relax 
\theauthors\else\theshortauthors\fi\fi\fi}}\def\@evenhead{\@oddhead}
\def\@oddfoot{\small\lfoot\ifnum\count0=\startpage\copyright\ \gtp\hfill\else
\agt, Volume \thevolumenumber\ (\thevolumeyear)\hfill\fi}
\def\@evenfoot{\@oddfoot}
\makeatother
\fi
\let\maketitlepage\makeagttitle
\let\makeshorttitle\maketitlepage
\let\maketitle\maketitlepage

\newwrite\gtoutfile
\long\gdef\makeheadfile{  
{\def\\{, }\def\s{ }
\immediate\openout\gtoutfile head.xxx
\immediate\write\gtoutfile{To: math@arxiv.org}
\immediate\write\gtoutfile{Subject: put OR rep NNNNN:ppppp}
\immediate\write\gtoutfile{--text follows this line--}
\immediate\write\gtoutfile{Proxy-for: \ifx\theasciiauthors\relax
\theauthors\else\theasciiauthors\fi\s<\ifx\theasciiemail\relax\theemail\else\theasciiemail\fi>}
\immediate\write\gtoutfile{\noexpand\\}
\immediate\write\gtoutfile{Authors: \ifx\theasciiauthors\relax
\theauthors\else\theasciiauthors\fi}
{\def\\{ }\immediate\write\gtoutfile{Title: \ifx\theasciititle\relax
\thetitle\else\theasciititle\fi}}
\immediate\write\gtoutfile{Subj-class: GT or SG, GR etc}
\immediate\write\gtoutfile{MSC-class: \theprimaryclass\ifx\thesecondaryclass\relax\else, \thesecondaryclass\fi}
\immediate\write\gtoutfile{Journal-ref: Algebr. Geom. Topol. \thevolumenumber\s
(\thevolumeyear) \startpage-\finishpage}
\immediate\write\gtoutfile{Comments: Published by Algebraic and
Geometric Topology at}
\immediate\write\gtoutfile{\s\s\s  http://www.maths.warwick.ac.uk/agt/AGTVol\thevolumenumber/agt-\thevolumenumber-\thepapernumber.abs.html}
\immediate\write\gtoutfile{\noexpand\\}
\immediate\write\gtoutfile{}
\ifx\theasciiabstract\relax
\immediate\write\gtoutfile{\theabstract}\else
\immediate\write\gtoutfile{\theasciiabstract}\fi
\immediate\write\gtoutfile{}
\immediate\write\gtoutfile{\noexpand\\}
\immediate\write\gtoutfile{}
\immediate\closeout\gtoutfile}}  

\def\maketitlepage{\makeagttitle\makeheadfile}
\let\makeshorttitle\maketitlepage
\let\maketitle\maketitlepage

\lognumber{22}
\volumenumber{3}
\volumeyear{2003}
\papernumber{22}
\published{25 June 2003}
\pagenumbers{623}{675}
\received{30 January 2003}
\accepted{4 June 2003}

\usepackage{amsmath,amssymb, pb-diagram, Bdiagx, mathrsfs}

\newcommand{\al}{\alpha}
\newcommand{\be}{\beta}
\newcommand{\ga}{\gamma}

\newcommand{\ep}{\epsilon}

\newcommand{\la}{\lambda}
\newcommand{\om}{\omega}

\newcommand{\Om}{\Omega}
\newcommand{\Si}{\Sigma}
\newcommand{\etab}{\tilde{\eta}}
\newcommand{\tsig}{\tilde{\sigma}}
\newcommand{\hal}{\hat{\alpha}}
\newcommand\SST{\scriptstyle}

\newcommand\tr{\mathop {\rm tr}}

\newcommand{\zz}{{\mathbb Z}}
\newcommand{\rr}{{\mathbb R}}
\newcommand{\cc}{{\mathbb C}}
\newcommand{\qq}{{\mathbb Q}}

\newcommand{\cE}{{\mathscr E}}
\newcommand{\cD}{{\mathscr D}}

\newcommand{\cL}{{\mathscr L}}

\newcommand{\cF}{{\mathscr F}}
\newcommand{\cG}{{\mathscr G}}
\newcommand{\cM}{{\mathscr M}}
\newcommand{\cH}{{\mathscr H}}
\newcommand{\cT}{{\mathscr T}}

\newcommand{\cP}{{\mathscr P}}

\newcommand{\Spec}{\operatorname{Spec}}
\newcommand{\Hom}{\operatorname{Hom}}

\newcommand{\image}{\operatorname{image}}
\newcommand{\End}{\operatorname{End}}
\newcommand{\Sign}{\operatorname{Sign}}
\newcommand{\SF}{\operatorname{SF}}
\def\eqref#1{(\ref{#1})}

\newcommand{\del}{\partial}

\newcommand{\chern}{\operatorname{ch}}
\renewcommand{\eqref}[1]{\textnormal{(\ref{#1})}}

\theoremstyle{definition}
\newtheorem{defn}{Definition}[section]

\theoremstyle{plain}
\newtheorem{lemma}[defn]{Lemma}
\newtheorem{thm}[defn]{Theorem}
\newtheorem{prop}[defn]{Proposition}
   \newtheorem{cor}[defn]{Corollary}

\newtheorem*{thmast}{Theorem}

\newcommand{\exendproof}{\renewcommand{\qed}{\relax}\end{proof}}

\begin{document}
\title{On the rho invariant for manifolds with boundary}

\authors{Paul Kirk\\Matthias Lesch}

\address{Department of Mathematics, Indiana University\\
Bloomington, IN 47405, USA}
\secondaddress{Universit\"at zu K\"oln, Mathematisches Institut\\
Weyertal 86--90, 50931 K\"oln, Germany}

\asciiaddress{Department of Mathematics, Indiana University\\
Bloomington, IN 47405, USA\\and\\Universitat zu Koln, 
Mathematisches Institut\\
Weyertal 86-90, 50931 Koln, Germany}

\email{pkirk@indiana.edu, lesch@mi.uni-koeln.de}

\url{http://php.indiana.edu/\char'176pkirk, 
http://www.mi.uni-koeln.de/\char'176lesch}

\begin{abstract}This article is a follow up of the previous article
 of the authors on the analytic surgery of $\eta$-- and
$\rho$--invariants. We investigate in detail the
(Atiyah--Patodi--Singer)--$\rho$--invariant for manifolds with
boundary.  First we generalize the cut--and--paste formula to
arbitrary boundary conditions. A priori the $\rho$--invariant is an
invariant of the Riemannian structure and a representation of the
fundamental group. We show, however, that the dependence on the metric
is only very mild: it is independent of the metric in the interior and
the dependence on the metric on the boundary is only up to its
pseudo--isotopy class. Furthermore, we show that this cannot be
improved: we give explicit examples and a theoretical argument that
different metrics on the boundary in general give rise to different
$\rho$--invariants. Theoretically, this follows from an interpretation
of the exponentiated $\rho$--invariant as a covariantly constant
section of a determinant bundle over a certain moduli space of flat
connections and Riemannian metrics on the boundary.  Finally we extend
to manifolds with boundary the results of Farber--Levine--Weinberger
concerning the homotopy invariance of the $\rho$--invariant and
spectral flow of the odd signature operator.
\end{abstract}

\asciiabstract{This article is a follow up of the previous article of
the authors on the analytic surgery of eta- and
rho-invariants. We investigate in detail the
(Atiyah-Patodi-Singer)-rho-invariant for manifolds with
boundary.  First we generalize the cut-and-paste formula to
arbitrary boundary conditions. A priori the rho-invariant is an
invariant of the Riemannian structure and a representation of the
fundamental group. We show, however, that the dependence on the metric
is only very mild: it is independent of the metric in the interior and
the dependence on the metric on the boundary is only up to its
pseudo--isotopy class. Furthermore, we show that this cannot be
improved: we give explicit examples and a theoretical argument that
different metrics on the boundary in general give rise to different
rho-invariants. Theoretically, this follows from an interpretation
of the exponentiated rho-invariant as a covariantly constant
section of a determinant bundle over a certain moduli space of flat
connections and Riemannian metrics on the boundary.  Finally we extend
to manifolds with boundary the results of Farber-Levine-Weinberger
concerning the homotopy invariance of the rho-invariant and
spectral flow of the odd signature operator.}



\maketitle

\numberwithin{equation}{section}

\section{Introduction}
The $\rho$--invariant of a closed odd-dimensional manifold was
defined in \cite{APS-I} as a difference of two spectral invariants.
To a closed Riemannian manifold $M$ and a unitary representation
of its fundamental group $\al\co\pi_1(M)\to U(n)$ Atiyah, Patodi,
and Singer assigned the real number
$$\rho(M,\al)= \eta(D_\al,M)- \eta(D_\tau,M)$$
where $D_\al$ denotes the odd signature operator with coefficients
in the flat bundle determined by $\al$, and $D_\tau$ is similar
with respect to the trivial representation $\tau$, and $\eta(D)$ denotes the
regularized signature of a self-adjoint Dirac operator $D$, introduced in
\cite{APS-I}.  As a consequence of their index theorem they showed that
\begin{enumerate}
\item $\rho(M,\al)$ is independent of the choice of Riemannian
metric on $M$.
\item $\rho(M,\al)$ extends the signature defect, that is, if
$(M,\al)=\partial (W,\beta)$ for some manifold $W$ with unitary
representation $\beta\co\pi_1(W)\to U(n)$, then
$$\rho(M,\al)=n\Sign(W)-\Sign_\beta(W)$$
where $\Sign_\beta(W)$ denotes the signature of $W$ with local
   coefficients in the flat bundle determined by $\beta$.
\end{enumerate}

Thus $\rho(M,\al)$ is a fundamental smooth invariant, but it remains largely
mysterious since its general definition depends on the spectra of differential
operators on $M$.

In \cite{Kirk-Lesch} we defined the $\rho$--invariant
$\rho(M,\al,g)$ in the case when the boundary of $M$ is non-empty
and proved a  non-additivity formula as a consequence of our
cut-and-paste formula for $\eta$--invariants of Dirac operators: if $M$ is a
closed manifold split into two parts $X$ and $Y$ along a hypersurface,
$$\rho(M,\al)=\rho(X,\al,g) + \rho(Y,\al,g) + m(V_{X,\al},V_{Y,\al})_{(g,\al)}
-m(V_{X,\tau},V_{Y,\tau})_{(g,\tau)}.$$
(For the definitions of the terms see \eqref{eq1.4} and \eqref{defofm}.)

   It is
the purpose of this article to explore the properties of
$\rho(M,\al,g)$, particularly those which flow from this formula. We will
describe the behavior of this invariant with respect to variations  in
Atiyah-Patodi-Singer (APS) boundary conditions, bordisms, and variations of
$\alpha$ and
$g$. As applications we prove   generalizations of the main results of
Farber-Levine-Weinberger \cite{Farber-Levine} concerning the homotopy
invariance of $\rho$ to manifolds with boundary. Special attention is given to
the construction of explicit examples.

The invariant $\rho(M,\al,g)$
is   defined as a difference of $\eta$--invariants for manifolds with boundary
and as such is also a spectral invariant. This has the   happy consequence
that it is gauge and isometry invariant.
But in contrast to the closed case,
when the boundary of $M$ is non-empty the resulting invariant
depends on the choice of Riemannian metric, $g$,  on the boundary.

Hidden
from the notation is the fact that elliptic boundary conditions
are required to define $\eta$--invariants on manifolds with
boundary. Our choice in \cite{Kirk-Lesch} is to use APS boundary
conditions with respect to the  Lagrangian subspace of {\em limiting
values of extended $L^2$ solutions} in the sense of \cite{APS-I}. This
   choice is intrinsic, homotopy invariant, and natural in a sense we
   will describe with respect to
   bordisms,  but is not continuous in families. This fact is apparent
   when one considers families for which the dimension of the kernel of the
   tangential operator is not constant, but discontinuities can also occur in
families
   for which the kernel is constant dimensional.

   For this reason it turns out to be useful to  allow more general
Lagrangian subspaces; we describe this generalization and derive
the corresponding cut-and-paste formula  for $\eta$- and $\rho$--invariants
with respect to arbitrary APS  boundary conditions
in Theorem \ref{thm1.1}.   Among other things, Theorem \ref{thm1.1} says:

\begin{thmast}
Suppose that $M= Y\cup_\Si X,$
   $\al\co\pi_1(M)\to U(n)$ is a representation,
$W^X,W^Y\subset H^*(\Si;\cc^n_\al)$ are  Lagrangian subspaces,
   and let $B$ be a flat connection on $M$ with holonomy $\al$ in
cylindrical form
near $\Si$. Then  the difference $\eta(D_B,M)-\eta(D_{B,W^X},X)
-\eta(D_{B, W^Y},Y)-m(W^X , W^Y )_{(\al,g)}$
equals the integer
\[
\tsig(V_{X,\al},V_{Y,\al},\ga(W^Y ))-\tsig(\ga(V_{X,\al}),W^X , W^Y ).
\]
\end{thmast}

In this statement $\eta(D_{B,W},X)$ denotes the $\eta$--invariant
of the odd signature operator coupled to a flat connection $B$   on
the manifold
$X$ with respect to APS boundary conditions determined by the Lagrangian
subspace
$W$ of the kernel of the tangential operator.  Moreover, $m(V,W)_{(\al,g)}$ is
an explicit real-valued invariant of pairs of Lagrangian subspaces of
the Hermitian
symplectic space
$H^*(\Si;\cc^n_\al)$ with its induced $L^2$ metric (this is defined in Section
\ref{Section2}), $\ga$ is the associated complex structure,
$\tsig$ is the Maslov triple index which appears in Wall's non-additivity
theorem
\cite{Wall},
$V_{X,\al}=\text{image}\big(H^*(X;\cc^n_\al)\to  H^*(\Sigma;\cc^n_\al)\big)$,
and similarly for $V_{Y,\al}$.

   This theorem   generalizes
\cite[Theorem 8.8]{Kirk-Lesch} to arbitrary APS boundary conditions.
Taking differences gives a corresponding formula for $\rho$--invariants.

   We next
give a topological description of how the spaces $V_{X,\al}$ propagate across a
bordism; the result is given as Theorem \ref{lem2.1} which gives a functorial
framework to keep track of APS boundary conditions and a companion additivity
formula for $\eta$ and $\rho$.

With these technical results in place, we can then begin a careful
investigation of how the $\rho$--invariants for manifolds with
boundary depend on
the choice of metric on the boundary and the representation.  For example, we
show:
\newtheorem*{Corfivetwo}{Corollary \ref{pseudo}}

\begin{Corfivetwo}The $\rho$--invariant for a
manifold with boundary depends on the Riemannian metric on the boundary only
up to its pseudo-isotopy class. Precisely,
if $f_0, f_1\co\del X\to \del X$ are  pseudo-isotopic  diffeomorphisms, then
$$
\rho(X,\al,f^*_0(g))=\rho(X,\al,f^*_1(g)).$$ 
\end{Corfivetwo}

In Section \ref{Dependence on the metric} we give explicit examples which show
that $\rho(X,\al,g)$ and $m(V,W)$ depend on the choice of Riemannian metric:

\newtheorem*{thmsixone}{Theorem \ref{dependence}}

\begin{thmsixone}\indent\par
\begin{enumerate}
\item There exists a 3-manifold $Y$ with non-empty boundary,  Riemannian
metrics $g_0$, $g_1$ on $\del Y$, and a representation $\al\co\pi_1Y\to
U(2)$ so that
$$\rho(Y,\al,g_0)\ne \rho(Y,\al ,g_1).$$ 
Examples exist with vanishing kernel of the tangential operator, i.e.\  
\[\ker A_b\cong H^*(\del Y;\cc^2_\al)=0.\]
\item There exist metrics $g_0$ and $g_1$ on the torus $T$ and 3-manifolds
$X$ and $Y$ with boundary $T$  such that setting
$V_X=  \text{image}\big( H^*(X;\cc)\to H^*(T;\cc)\big)$ and
$V_Y= \text{image}\big(  H^*(Y;\cc)\to H^*(T;\cc)\big)$ (with
$\theta$ the trivial conection),
$$m(V_X,V_Y)_{(\theta, g_0)}\ne m(V_X,V_Y)_{(\theta, g_1)}.$$
\end{enumerate}
\end{thmsixone}


This theorem drives home the point that the choice of Riemannian metric on the
boundary is an essential ingredient of the $\rho$--invariant on a manifold with
non-empty boundary.

In Section \ref{flwtheorem} we extend to
manifolds with boundary the results of
Farber-Levine-Weinberger  concerning the homotopy invariance of the
$\rho$--invariants and spectral flow of the odd signature operator.
Let $$\chi(\pi_1 X ,
U(n))=\Hom(\pi_1 X ,U(n))/\text{conjugation} $$ and let $\mathscr{M}_{\del X}$
denote the space of Riemannian metrics on $\del X$.  Notice that a map $F\co X\to
X'$ which restricts to a diffeomorphism on the boundary and which induces an
isomorphism on fundamental groups provides an identification of
$\mathscr{M}_{\del X}$ with $\mathscr{M}_{\del X'}$ and $\chi(\pi_1 X ,U(n))$
with
$\chi(\pi_1X',U(n))$.

\newtheorem*{thmseventwo}{Theorem \ref{FLWKL}}

\begin{thmseventwo}Let $F\co X\to X'$ be a homotopy
equivalence of compact manifolds which restricts to a diffeomorphism on the
boundary.   Then the difference
   $$\rho(X)-\rho(X')\co\chi(\pi_1X,U(n))\times \mathscr{M}_{\del X}\to \rr$$
factors through  $\pi_0(\chi(\pi_1X, U(n)))\times (\mathscr{M}_{\del X}/
\mathscr{D}^0_{\del X})$
(where $\mathscr{D}^0_{\del X}$ denotes the group of diffeomorphisms
of $\del X$ pseudo-isotopic to the identity)  and takes values in
the rational numbers.

In other words there is a commutative diagram
\[\begin{diagram} \dgmag{1000}\dgsquash
\node{\chi(\pi_1X,U(n)) \times \mathscr{M}_{\del X}}\arrow{e,t}{\rho(X)-
\rho(X')}\arrow{s}\node{\rr}\\
\node{\pi_0(\chi(\pi_1X,U(n))) \times (\mathscr{M}_{\del X}/
\mathscr{D}^0_{\del X})}\arrow{e}\node{\qq}\arrow{n}
\end{diagram}\]
Moreover the difference $\rho(X, \al,g)-\rho(X',\al, g)$ vanishes for $\al$ in
the path component of the  trivial representation.
\end{thmseventwo}

We also show that the spectral flow of the odd signature operator coupled to
a   path of flat connections is a homotopy invariant. In the following
statement $P(t)$ denotes a smooth path of self-adjoint APS boundary
conditions with prescribed endpoints. (We show how to construct such a path in
Lemma
\ref{apslem}.)
\newtheorem*{thmsevenfour}{Theorem \ref{sfthm}}

\begin{thmsevenfour}Suppose that $F\co X'\to X$ is
a homotopy equivalence
which restricts to a diffeomorphism $f=F|_{\del X'}\co\del X'\to \del X$.  Assume
that $B_t$ is a continuous, piecewise smooth path of flat $U(n)$  connections
on $E\to X$. Use $F$ to pull back the path $B_t$ to a path of flat  connections
$B'_t$ on $X'$ and to identify $\del X$ with $\del X'$, and choose a path
$P(t)$ of APS boundary conditions as in Lemma \ref{apslem}.

Then
$$\SF(X,D_{B_t,P(t)})_{t\in [0,1]}=  \SF(X',D_{B'_t,P(t)})_{t\in [0,1]}.$$
\end{thmsevenfour}


   In Section \ref{det bundles} we   use the machinery
of determinant bundles, especially the Dai-Freed theorem, to study the
variation of $\rho(X,\al,g)$ modulo $\zz$.  By working modulo $\zz$ one loses
   geometric information but   the discontinuities of
$\rho$ as a function of $\al$ are eliminated.  In particular variational
techniques can be applied.

      Theorem
\ref{thmcov} implies the following. In this theorem $\nabla^Q$
denotes the connection on the determinant bundle as introduced by
Quillen in \cite{Quillen}. (See \cite{BGV} for the
construction of $\nabla^Q$ in general.)

\begin{thmast}The assignment of the exponentiated $\rho$--invariant to a flat
$SU(n)$ connection
$B$ on a manifold with boundary $X$ and a choice of Riemannian metric $g$ on
$\del X$,
$$(B,g)\mapsto \exp(\pi i\
\rho(X,\al,g)),$$
(where $\al$ is the holonomy of $B$) defines a smooth horizontal (with respect to the   connection
$\nabla^Q$)
cross section   of the
determinant bundle of the family of tangential operators to the odd signature
operators.  
\end{thmast}

This theorem allows one to relate the  mod $\zz$ reduction of the
$\rho$--invariant on manifolds with isomorphic fundamental groups and
diffeomorphic boundaries, and also shows that the manner in which
$\rho(X,\al,g)$ depends on the choice of metric $g$ on $\del X$ is intimately
tied to the
  connection $\nabla^Q$.

For example, the following is a consequence of Theorem \ref{cor16}.  We view
$\rho(X)$ as a function of the conjugacy class of the representation $\al$ and
the metric $g$.

\begin{thmast}Let $X$ and $X'$ be two odd dimensional
manifolds and suppose that $F\co X'\to X$ is a smooth map which induces
an isomorphism on fundamental groups  and such that the restriction
$f=F|_{\del X'}\co\del X'\to \del X$ is a diffeomorphism.

Then  there is a factorization
\[\begin{diagram}\dgmag{600}
\node{\chi(\pi_1(X), SU(n))\times
\cM_\Sigma }\arrow[2]{e,t}{\rho(X)-\rho(X')}\arrow{se}\node[2]{\rr/\zz}
\\ \node[2]{\pi_0(\chi(\pi_1(X),
SU(n)))\times\cM_\Sigma }\arrow{ne,..}
\end{diagram}\]
and $\rho(X)-\rho(X')$ is zero on the path component of the
trivial representation. The result holds for $U(n)$ replacing
$SU(n)$ if $\dim X=4\ell-1$.
\end{thmast}

These results, together with the cut--and--paste formula for $\rho$--invariants
(Theorem \ref{thm1.1}) are a step in the program of determining what the
homotopy
    properties of the $\rho$--invariant are. A discussion of problems in
this topic
is given in Section \ref{Section 9}, including the following consequence of
Theorem \ref{split} concerning the homotopy invariance of the $\rho$--invariant
for closed manifolds.


\begin{thmast} Let $M$ and $M'$ be closed
manifolds, and suppose there exists a sep\-a\-rating hypersurface $\Si\subset M$
and a smooth homotopy equivalence
$F\co M'\to M$ so that the restriction of $F$ to $F^{-1}(\Si)$ is a
diffeomorphism.  Write $M=X\cup_\Si Y$ and $M'=X'\cup_\Si Y'$ and suppose that
$F$ restricts to homotopy equivalences $X'\to X$ and $Y'\to Y$.  Let $\al\in
\chi(\pi_1M, U(n))$.

If the restriction  $\al|_{X}$ (resp. $\al|_Y$) of $\al$ to
$\pi_1X$ (resp. $\pi_1 Y$) lies in the path component of the
trivial representation of $\chi(\pi_1X,U(n))$ (resp.
$\chi(\pi_1Y, U(n))$)    then $\rho(M,\al)=\rho(M',\al)$.
\end{thmast}

We finish the article with a brief discussion of the relation of our
investigations to one of the approaches to the program of constructing
topological quantum field theories proposed in \cite{atiyah-GPK}.

\rk{Acknowledgements}

The   authors  thank J.F.~Davis,   C.~Livingston,
S.~Paycha, and K.P.~Wojciechowski for very helpful discussions.

The first named author gratefully acknowledges the support of
the National Science Foundation under grant no. DMS-0202148.

\section{The $\rho$--invariant on manifolds with boundary}\label{Section2} We begin
by recalling the context and the definition of the $\rho$--invariant  for a
manifold with boundary. More details can be found in \cite{Kirk-Lesch}.

Let $X$ be a $(2k+1)$-dimensional smooth, oriented, compact manifold with
(possibly empty) boundary $\Si$.  Fix a Riemannian metric $\tilde{g}$ on
$X$  in product form near the boundary $\Si$.  To keep track
of signs it is crucial to fix a convention for the orientation
of a collar of the boundary. In this paper we will use the convention
of \cite{Kirk-Lesch}:\label{signstuff} if not indicated otherwise a collar of the boundary
will be written as $\Sigma\times [0,\ep)$, i.e.\ the manifold $X$ is
``on the right'' of the boundary.   The choice of
the sign convention has consequences for the definition of $A_b$ and $\gamma$
(and hence
the Hermitian symplectic structure on $H^*(\Sigma;\cc^n_\al)$) below.

Let $B$ be a flat $U(n)$
connection on
$X$ in product form near the boundary, i.e.\ $B|_{\Si\times
[0,\ep)}=\pi^*(b)$ for some flat connection  $b$ on $\Si$; here
$\pi\co\Si\times[0,\ep)\to \Si$ denotes the projection.   Denote by
$\al\co\pi_1X\to U(n)$ the holonomy representation of $B$.  Since it will
be central in what follows, denote the restriction of $\tilde{g}$ to the
boundary $\Si$ by $g$.

\begin{sloppy}
The {\em odd signature operator coupled to the flat connection $B$}
$$D_B\co\oplus_p\Om^{2p}(X;E)\to\oplus_p\Om^{2p}(X;E)$$ is defined by
$$D_B(\beta)=i^{k+1}(-1)^{p-1}(*d_B-d_B*)(\beta) \hbox{ for } \beta\in
\Om^{2p}(X;E).$$ 
Here,  $*\co\Om^{\ell}(X;E)\to
\Om^{2k+1-\ell}(X;E)$ denotes the Hodge * operator (which is determined by
the Riemannian metric
$\tilde{g}$ on $X$), $d_B\co\Om^{\ell}(X;E)\to
\Om^{\ell+1}(X;E)$ denotes the covariant derivative associated to the flat
connection $B$, and $E\to X$ denotes the associated
Hermitian $\cc^n$ vector bundle.

\end{sloppy}

On the collar $\Si\times [0,\ep)$, $D_B$ takes the form (after conjugating
with a certain unitary transformation, see \cite[(8.1)]{Kirk-Lesch} for details)
$$D_B= \gamma(\frac{\del}{\del x}  + A_b),$$ where the {\em de~Rham
operator}
$$A_b\co\oplus_k\Omega^{k}(\Si;E|_{\Si})\to
\oplus_k\Omega^{k}({\Si};E|_{{\Si}})$$
   is defined by $$A_b(\beta)=
\begin{cases} -(d_b\hat{*}+\hat{*}d_b)\beta, & \text{if }
\beta\in
\oplus_k\Omega^{2k}({\Si};E|_{{\Si}}),\\
\quad (d_b\hat{*}+\hat{*}d_b)\beta, & \text{if } \beta\in
\oplus_k\Omega^{2k+1}({\Si};E|_{{\Si}}).\end{cases}$$ In these formulas
$\hat{*}$ denotes the Hodge $*$ operator on
$\Si$ and $$\gamma\co\oplus_{p}\Omega^{p}({\Si }; E|_{{\Si}})\to
\oplus_{p}\Omega^{p}({\Si}; E|_{{\Si}})$$
    coincides with $\hat{*}$ up to a constant:
$$ \gamma(\beta)=\begin{cases} i^{k+1}(-1)^{p-1}\hat{*}\ \beta,  &
\text{if }
\beta\in
\Omega^{2p}({\Si}; E|_{{\Si}}), \\ i^{k+1}(-1)^{k-q}\hat{*}\ \beta, &
\text{if }
\beta\in
\Omega^{2q+1}({\Si}; E|_{{\Si}}). \end{cases}$$
One calculates that $\gamma^2=-Id, \gamma A_b=-A_b\gamma$,
   and that $\gamma$ is unitary with respect
to the $L^2$ inner product on
$\Om^*(\Si; E|_{\Si})$ defined by
$$\langle \beta_1,\beta_2\rangle=\int_\Si \beta_1\wedge \hat{*}\beta_2.$$ (The
   Riemannian metric on $\Si$ is used to define the Hodge $*$-operator
$\hat{*}$,
and we have suppressed the notation for the inner product in the bundle
$E$.) The operator $A_b$ is elliptic and self-adjoint and hence
one has an orthogonal decomposition
\begin{equation}\label{eq1.2}L^2(\Om^*(\Si;E|_{\Si}))= F_b^-\oplus \ker
A_b\oplus F_b^+\end{equation} into the negative eigenspan, kernel, and
positive eigenspan of
$A_b$. The relation $\gamma A_b=-A_b\gamma$ implies that
$\ker A_b$ is preserved by $\gamma$ and that
$\gamma$ maps  $F_b^+$ unitarily onto $F_b^-$.


   The kernel of  $A_b$ is identified by the Hodge theorem with the twisted
de~Rham  cohomology of the complex $(\Om^*(\Si; E|_{\Si}),d_b)$;
indeed the elements of $\ker A_b$ are just the $d_b$-harmonic
forms and so the composite $$\ker A_b=\ker d_b\cap \ker
d_b^*\subset \ker d_b\to \frac{\ker d_b}{\image d_b}$$ is an
isomorphism. The de~Rham theorem then identifies the
cohomology of $(\Om^*(\Si; E|_{\Si}),d_b)$ and the (singular or cellular)
cohomology
$H^*(\Si;\cc^n_\al)$ with local coefficients given by the holonomy
representation $\al$.

The triple $(\ker A_b, \langle\ , \ \rangle, \ga)$ gives $\ker A_b$ the
structure of  a  Hermitian symplectic space. In general
a {\em Hermitian symplectic space} $(H,\langle\ , \ \rangle, \ga)$
is a finite dimensional complex vector space $H$ with a positive
definite Hermitian inner product $\langle\ , \ \rangle \co H\times H\to \cc$
and an isomorphism $\ga\co H\to H$ which is unitary, i.e.\
$\langle\ga(x),\ga(y) \rangle
=\langle x , y \rangle$,  satisfying $\ga^2=-I$ such that the signature of
$i\ga$ is zero.   The {\em underlying symplectic structure}
   is the pair $(H,\omega)$, where $\omega$ is the
   non-degenerate skew-Hermitian form
$$\omega(x,y)=\langle x, \ga(y)\rangle.$$
The signature of $i \ga$ on $\ker A_b\cong
H^*(\Si;\cc^n_\al)$ is zero. This is a
consequence of the fact that $(\Si,\al|_\Si)$ bounds $(X,\al)$, and is
not true for a general pair $(\Si,\al)$.  However, it is true in many
important cases, for example if $\Si$ is a $4\ell-2$ dimensional manifold and
$\al\co\pi_1\Si\to U(n)$ factors through $O(n)$.

   In contrast to the Hermitian inner
product $\langle\ , \ \rangle$ and the unitary map $\ga$ on $\ker A_b$,
the symplectic form $\omega$ does not depend on the
Riemannian metric and in fact is given by the cup product:
$$\om(\beta_1,\beta_2)=i^r \int_\Si \beta_1\wedge
\beta_2=i^r([\beta_1]\cup[\beta_2])\cap [\Si], $$ where $i^r$ is a
constant depending on the degrees of the $\beta_i$.

A subspace $W$ of a Hermitian symplectic space $(H,\om)$ is called
{\em Lagrangian} if $\omega$ vanishes on $W$  and $W$ is maximal
with this property.
    This is equivalent to $\ga(W)=W^\perp$, but being a Lagrangian subspace is a
property of the underlying symplectic structure.
Note that $\dim W=\frac{1}{2}\dim H$.
Denote the Grassmannian of all Lagrangian subspaces of $H$ by
$\cL(H)$.

We summarize:  The symplectic structure on $H^*(\Si;\cc^n_\al)$,
and hence the Grassmannian $\cL(H^*(\Si;\cc^n_\al))$, depends only
on the cohomology and cup product, and therefore is a homotopy
invariant of $(\Si,\al)$. On the other hand, the Hermitian
symplectic structure on $H^*(\Si;\cc^n_\al)$ depends on its
identification with $\ker A_b$ via the Hodge and de Rham theorems,
since the inner product $\langle \ , \ \rangle$ is the restriction
of the $L^2$ inner product (which depends on the Riemannian metric
on $\Si$) to $\ker A_b$.

The following lemma is well--known; it follows by a standard argument
using Poincar{\'e} duality (cf.\ also \cite[Cor.\ 8.4]{Kirk-Lesch}).
\begin{lemma}\label{lemma1.2} The image of the restriction
\begin{equation}\label{eq1.3}H^*(X;\cc^n_\al)\to H^*(\Si; \cc^n_\al)
\end{equation} is a Lagrangian subspace.\qed
\end{lemma} We will denote this subspace by $V_{X,\al}$, and, by slight
abuse of notation,  its preimage in $\ker A_b$ via the isomorphism
$\ker A_b\cong H^*(\Si;\cc^n_\al)$ will also be denoted by $V_{X,\al}$. We
emphasize  that the
Lagrangian
$V_{X,\al}$ is a homotopy  invariant of $(X,\al)$. Moreover it gives a
distinguished
element in the Grassmannian  $ \mathscr{L}(H^*(\Si;
\cc^n_\al))$. Considered as a subspace of  $\ker A_b$,
$V_{X,\al}$ coincides with  the {\em limiting
values of extended $L^2$ solutions  of} $D_B\phi=0$ on $(\Si\times(-\infty,0])\cup X$
 in the sense of \cite{APS-I}.

Lagrangian subspaces  of  $ H^*(\Si; \cc^n_\al)$ are used to
produce elliptic self-adjoint Atiyah-Patodi-Singer (APS) boundary
conditions for the odd signature operator $D_B$ as follows. Given a
Lagrangian subspace $W\subset H^*(\Si; \cc^n_\al)$ we consider the
orthogonal projection in $L^2(\Om^*(\Si;E|_{\Si}))$ onto $F_b^+\oplus
W$. This orthogonal projection defines a well--posed boundary
condition for $D_B$ (see e.g.~\cite{BooWoj}).

Restricting $D_B$ to the subspace of sections whose restriction to the
boundary lies in the kernel of this projection makes $D_B$ a discrete
self-adjoint operator which we denote by $D_{B,W}$.   The following properties
of this operator are the starting point of the investigations of
this article and go back to Atiyah, Patodi, and Singer's
fundamental articles \cite{APS-I, APS-II, APS-III}.  In this
context the following facts are explained in \cite{Kirk-Lesch}.

\begin{enumerate}
\item The $\eta$ function of the operator $D_{B,W}$,
$$\eta(s)=\sum_{\la\in \text{\tiny Spec}(D_{B,W})\setminus\{0\}}
\operatorname{sign}{\la}\,|\la|^{-s},$$ converges for Re$(s)>>0$ and has a
meromorphic continuation to the entire complex plane with no pole
at $s=0$. Denote its value at $s=0$ by $$\eta(D_{B,W},X):=\eta(0).$$
\item  The kernel  of  $D_{B,W}$   fits into an exact sequence
\begin{equation}\label{eq1.6}\begin{split}
0\to \big(\image H^*(X,\Si; \cc^n_\al)&\to H^*( X; \cc^n_\al)\big)
\to \ldots\\ \ldots&\to\ker D_{B,W}\to W\cap \ga(V_{X,\al})\to 0.
			     \end{split}
\end{equation}
In particular, taking $W=V_{X,\al}$ we see $$\ker D_{B,V_{X,\al}}=
\image H^*(X,\Si; \cc^n_\al)\to H^*( X; \cc^n_\al).$$
\end{enumerate}

We next recall the definition of the $\rho$--invariant
for manifolds with boundary from  \cite{Kirk-Lesch}.   Let $\Theta$
denote the trivial connection in the product bundle $\cc^n\times X$ in the
form $\Theta=\pi^*(\theta)$ in the collar of $\partial X$,  and
$\tau\co\pi_1X\to U(n)$ the trivial representation.  Then define
\begin{equation}\label{eq1.4}
\rho(X,\al,g)=\eta(D_{B,V_{X,\al}},X)- \eta(D_{\Theta,V_{X,\tau}}, X).
\end{equation}
It is shown in \cite[Sec.\ 8]{Kirk-Lesch} that $\rho(X,\al,g)$
depends only on the smooth structure on $X$, the conjugacy class of
$\al\co\pi_1X\to U(n)$, and the Riemannian metric $g$ on $\Si=\del X$. In
particular, it is independent of the choice of flat connection $B$ with
holonomy conjugate to
$\al$ and also independent  of the Riemannian metric $\tilde{g}$ on $X$
extending $g$.

When $\del X$ is empty, then the diffeomorphism invariance of
$\rho(X,\al)$ was established by Atiyah, Patodi, and Singer in
\cite{APS-II} and follows straightforwardly from their index
theorem.   The cut--and--paste formulae
\begin{equation}\label{eq1.49}
\eta (D_B,M)
      =\eta(D_{B,V_{X,\al}},X)+\eta(D_{B,V_{Y,\al}},Y)
     +
m(V_{X,\al},
V_{Y,\al})_{(b,g)}
\end{equation}
and
\begin{equation}\label{eq1.5}
\rho (M,\al)
      =\rho(X,\al,g)+\rho(Y,\al,g)
     +
m(V_{X,\al},
V_{Y,\al})_{(b,g)}-m(V_{X,\tau},
V_{Y,\tau})_{(\theta,g)}
\end{equation}
when $M=Y\cup_\Si X$
were proven in \cite[Sec.\ 8]{Kirk-Lesch} and are
the basis for our investigations in the present article.

In Equations \eqref{eq1.49} and
\eqref{eq1.5} the correction term $ m(V,W)_{(b,g)}$ is  a real
valued invariant of pairs of Lagrangians in $H^*(\Si;\cc^n_\al)$;
it depends on the identification of
$H^*(\Si;\cc^n_\al)$ with the kernel of $A_b$ and hence may {\it a
priori} (and {\it a posteriori} as well, see Section \ref{Dependence on the
metric})
depend on the Riemannian metric $g$ on $\Si$. It is defined as follows.

\begin{sloppypar}
Let $\ker A_b^+$ denote the $+i$-eigenspace  of $\ga$
   acting on $\ker A_b$ and let
$\ker A_b^-$ denote the 
$-i$-eigenspace.  Then every Lagrangian subspace $W$ of
   $\ker A_b\cong H^*(\Si;\cc^n_\al)$ can be written uniquely as a
   graph
\begin{equation}\label{eqforphi}W=\{ x +
\phi(W)(x) | x\in \ker A_b^+\},\end{equation}
where $\phi(W)\co\ker A_b^+\to \ker A_b^-$ is a unitary isomorphism.
   The map $W\mapsto \phi(W)$
determines a diffeomorphism between the space $\mathscr{L}(\ker A_b)$ of
Lagrangians in
$\ker A_b$ to the space of unitary
operators $U(\ker A_b^+, \ker A_b^-)$.  We take the branch $\log(re^{it})=\ln
r + it,\  r
>0, -\pi<t\leq \pi$ and use it to define 
$\tr\log\co U(\ker A_b^+)\to i \rr$ via
   $\tr\log(U)=\sum \log(\lambda_i), \lambda_i\in \Spec U$. Then define
\begin{equation}\label{defofm}
\begin{split}
m(V,W)_{(b,g)}&=  -\tfrac{1}{\pi i}\tr\log (-
\phi(V)\phi(W)^*)
+\dim(V\cap W)\\
&=-\tfrac{1}{\pi i}
\sum_{\begin{array}{c}\SST\lambda\in\Spec(-\phi(V)\phi(W)^*)\\
\SST\lambda\neq -1\end{array}} \log \lambda.
  \end{split} \end{equation}  We will abbreviate this to $m(V,W)$
when $(b,g)$ is clear from context. Since $-\phi(V)\phi(W)^*$ is unitary, its eigenvalues
are unit complex numbers, and hence $m(V,W)$ is a real number. The term $\dim(V\cap W)$ is added to match conventions and to simplify formulas; notice that its effect is to remove the contribution of the $-1$ eigenspace of $-\phi(V)\phi(W)^*$ to $\tr\log (-\phi(V)\phi(W)^*)$. Thus $m$ is not in general a continuous function of $V$ and $W$.
The function $m$ has been investigated before, the notation is taken from
\cite{Bun:GPE}.
\end{sloppypar}

\section{Cutting and pasting formulas with arbitrary\newline boundary conditions}

\begin{sloppypar}
The $\eta$--invariants appearing in the definition of $\rho$ of Equation
   \eqref{eq1.4} are taken with respect to the boundary conditions
$V_{X,\al}\subset H^*(\Si;\cc^n_\al)$ and $V_{X,\tau}\subset H^*(\Si;\cc^n)$.
   More precisely, the Lagrangian $V_{X,\al}\subset H^*(\Si;\cc^n_\al)$
determines
a subspace (still denoted $V_{X,\al}$) of $\ker A_b$, and this in turn
determines the orthogonal projection to $F_b^+\oplus V_{X,\al}$, (recall that
$F_b^+$ is shorthand for the positive eigenspan of $A_b$). A similar comment
applies to
   $V_{X,\tau}$.  Since these Lagrangians are canonically
determined  by the homotopy type of the pair $(X,\al)$ and the Riemannian
metric on
$\Si$, they present a natural choice for the boundary conditions. Nevertheless
it is useful to use other Lagrangians in $H^*(\Si;\cc^n_\al)$ to define
boundary conditions. One important reason is that the $V_{X,\al}$ do not vary
continuously in families, even if $\ker A_b$ does.
\end{sloppypar}

\begin{defn}\label{generalrho} Let $X$ have boundary $\Si$ and let
$\al\co\pi_1X\to U(n)$ be a representation.  Given
   Lagrangian subspaces  $W_\al\subset H^*(\Si;\cc^n_\al)$ and
$W_\tau\subset H^*(\Si;\cc^n)$, define
$\rho(X,\al,g,W_\al,W_\tau)$ by
$$\rho(X,\al,g,W_\al,W_\tau):=\eta(D_{B,W_\al},X)-\eta(D_{\Theta,W_\tau},X).$$
\end{defn}
Thus $\rho(X,\al,g)$ is shorthand for $\rho(X,\al,g,V_{X,\al},V_{X,\tau})$.

We next recall the definition of $\tsig$ from
\cite[Sec.\ 8]{Kirk-Lesch}.  Given Lagrangian subspaces $U,V,W$ of a Hermitian
symplectic space, define
$$\tsig(U,V,W):=m(U,V) + m(V,W) + m(W,U).$$
Then $\tsig$ is  integer-valued, depends only on the symplectic form $\omega$,
and coincides with Wall's correction term
for the non-additivity of the signature \cite{Wall} as well as the Maslov
triple index of \cite{CLM}.

The following theorem gives a complete formulation of the dependence
of the $\eta$- and $\rho$--invariants for a manifold with boundary on the choice
of Lagrangians used for APS boundary conditions.

\begin{thm}\label{thm1.1} Suppose that $M= Y\cup_\Si X$,
   $\al\co\pi_1(M)\to U(n)$ is a representation,
$W^X_\al,W^Y_\al\subset H^*(\Si;\cc^n_\al)$ and $W^X_\tau,W^Y_\tau\subset
H^*(\Si;\cc^n)$ are  Lagrangian subspaces,
   and let $B$ be a flat connection on $M$ with holonomy $\al$ in
cylindrical form
near $\Si$. Orientation dependent quantities like $\gamma$ etc. are taken
with respect to $X$ according to the convention explained on page \pageref{signstuff}.

Then:
\begin{enumerate}
\item
$\eta(D_{B,W^X_\al},X)-\eta(D_{B,V_{X,\al}},X)= m(\gamma(V_{X,\al}),W^X_\al).$
\item $\rho(X,\al,g,W^X_\al,W^X_\tau)$
depends only on the diffeomorphism type of $X$, the representation
$\al$, the Lagrangian subspaces $W_\al^X,W_\tau^X$ and
   the Riemannian metric $g$ on $\Si=\del X$.
\item The difference $\eta(D_B,M)-\eta(D_{B,W^X_\al},X)
-\eta(D_{B,W^Y_\al},Y)-m(W^X_\al, W^Y_\al)$
is an integer. In fact it equals
\[
\tsig(V_{X,\al},V_{Y,\al},\ga(W^Y_\al))-\tsig(\ga(V_{X,\al}),W^X_\al, W^Y_\al).
\]
\item $\eta(D_B,M)= \eta(D_{B, V_{X,\al}},X) +
\eta(D_{B,\ga(V_{X,\al})},Y)$ and so
$$\rho(M,\al)= \rho(X,\al, V_{X,\al}, V_{X,\tau}) + \rho(Y,
\al,\ga(V_{X,\al}), \ga(V_{X,\tau})).$$
\end{enumerate}
\end{thm}

\begin{proof}We use the results of \cite{Kirk-Lesch}. Recall the notation
$\etab(D)=\tfrac{1}{2}(\eta(D) +\dim \ker D)$. For the proof of (1)
we omit the sub- and superscripts of $W$:

By \cite[Theorem
4.4]{Kirk-Lesch} we have
\begin{equation}\begin{split}
     \etab&(D_{B,W} ,X)-\etab(D_{B,V_{X,\al}},X)\\
        &=\tfrac{1}{2\pi i}
\bigl(\tr\log(\Phi(P^+(W))\Phi(P_X)^*)-
\tr\log(\Phi(P^+(V_{X,\al}))\Phi(P_X)^*)\bigr).
		\end{split}
\label{G2.7}
\end{equation}
Here, $P^+(W)$ denotes the orthogonal projection onto $W\oplus
F_b^+$, $P_X$ denotes the Calder\`on projector for $D_B$ acting on $X$,
and
$\Phi$ is the infinite--dimensional version of
$\phi$:
   it denotes the diffeomorphism from the
(infinite--dimensional) Lagrangian Grassmannian onto
$\mathscr{U}(\ker(\gamma-i),\ker(\gamma+i))$ (cf.\ \cite[Sec.\
2]{Kirk-Lesch}). Using \cite[Lemma 6.9]{Kirk-Lesch} we identify
the right side of \eqref{G2.7} with
\begin{equation}
    \tau_{\mu}(P^+(V_{X,\al}),P^+(W),P_X)-\tfrac{1}{2\pi i}\tr\log\bigl(
      \Phi(P^+(V_{X,\al}))\Phi(P^+(W))^*\bigr),\label{G2.8}
\end{equation}
where $\tau_\mu$ is the Maslov triple index defined in
      \cite[Sec.\ 6]{Kirk-Lesch}.

In view of \cite[Lemma 8.10]{Kirk-Lesch} the quantity
$\tau_{\mu}(P^+(V_{X,\al}),P^+(W),P_X)$ is invariant under
adiabatic stretching and equals
\begin{equation}
     \tau_\mu(V_{X,\al},W,V_{X,\al})=\dim\big( V_{X,\al}\cap\gamma (W)\big),
\end{equation}
where the last equality follows from \cite[Prop.\ 6.11]{Kirk-Lesch}.

As in the proof of \cite[Theorem 8.12]{Kirk-Lesch} one calculates
\begin{equation}\label{G2.9}
    \tr\log\bigl(\Phi(P^+(V_{X,\al}))\Phi(P^+(W)^*)\bigr)=
     \tr\log\bigl(\phi(V_{X,\al})\phi(W)^*\bigr).
\end{equation}
The identity $\ga^2=-I$ shows that $\dim(V_{X,\al}\cap
\ga(W))=\dim(\ga(V_{X,\al})\cap W)$ and clearly
$\phi(\ga(W))=-\phi(W)$.   These facts together with the
definition of $m(V,W)$  and Equation \eqref{G2.9} imply
\begin{equation}
\begin{split}\etab(D_{B,W},X)&-\etab(D_{B, V_{X,\al}},X)\\
&=\dim (\ga(V_{X,\al})\cap
W) -\tfrac{1}{2\pi i}\tr\log (\phi(V_{X,\al})\phi(W)^*)\\
&=\tfrac{1}{2}\big(m(\ga(V_{X,\al}),W) + \dim(\ga(V_{X,\al})\cap
W)\big).
\end{split}\end{equation}
Using the definition $\etab(D)= \tfrac{1}{2}(\eta(D)+\dim \ker D)$ we see that
$\eta(D_{B,W},X)-\eta(D_{B,V_{X,\al}}, X)- m(\gamma(V_{X,\al}),W)$ equals
\begin{equation}\label{G2.95}
-\dim \ker D_{B,W}+\dim \ker D_{B,V_{X,\al}} +\dim(\gamma(V_{X,\al})\cap W).
\end{equation}
But \eqref{G2.95} vanishes,
   as one sees by using the exact sequence \eqref{eq1.6}. This
   proves the first assertion of Theorem \ref{thm1.1}.

The second assertion follows from the first part
   and \cite[Lemma 8.15]{Kirk-Lesch}.

Using \eqref{eq1.49} and the first assertion one sees that
   \begin{equation*}\eta(D_B,M)-\eta(D_{B,W^X_\al},X)
-\eta(D_{B,W^Y_\al},Y)-m(W^X_\al, W^Y_\al)
\end{equation*}
equals
\begin{equation}\label{G2.91}
m(V_{X,\al}, V_{Y,\al}) - m(\ga(V_{X,\al}), W^X_\al) +
   m(\ga(V_{Y,\al}), W^Y_\al)-m(W^X_\al,W^Y_\al).\end{equation}
(There is one subtlety: the sign change of the term
$ m(\ga(V_{Y,\al}), W^Y_\al)$ occurs because viewed from the ``$Y$'' side, the
Hermitian symplectic structure changes sign.)

Using the identities $m(V,W)=-m(W,V)$ and $\phi(\ga(W))=-\phi(W)$, so that
$m(\ga(V),\ga(W))=m(V,W)$, we can rewrite \eqref{G2.91} as
\[
-m(\ga(V_{X,\al}),W^X_\al)-m(W^X_\al,W^Y_\al) + m(V_{X,\al}, V_{Y,\al})
+ m(V_{Y,\al}, \ga(W^Y_\al)),\]
which equals
\begin{equation*}\begin{split}
&\tsig(V_{X,\al},V_{Y,\al},\ga(W^Y_\al))-m(\ga(W^Y_\al),V_{X,\al})\\
&\quad-\tsig(\ga(V_{X,\al}), W^X_\al,W^Y_\al) + m(W^Y_\al,\ga(V_{X,\al}))\\
=&\tsig(V_{X,\al},V_{Y,\al},\ga(W^Y_\al))
-\tsig(\ga(V_{X,\al}), W^X_\al,W^Y_\al)
\end{split}
\end{equation*}
as desired. This proves the third assertion.

The last statement follows straightforwardly from the previous or,
alternatively, can be immediately recovered from
\cite[Theorem 8.8]{Kirk-Lesch}.
\end{proof}

\section{Lagrangians induced by bordisms}

Theorem \ref{thm1.1}   gives splitting formulas for the $\eta$ and
$\rho$--invariants of $D_B$ in the situation when a manifold $M$ is
decomposed into two pieces $X$ and $Y$ along a hypersurface
$\Si$.  To develop this into a useful cut-and-paste machinery for
the $\rho$--invariant requires keeping track of  the Lagrangian
subspaces $V_{X,\al}=\image H^*(X; \cc^n_\al)\to H^*(\Si;
\cc^n_\al)$ and their generalizations.  It is clearest to give an
exposition based on the effect of a bordism on Lagrangian subspaces
and we do this next.

   Let $X$ be a Riemannian manifold with boundary
   $-\Si_0\amalg\Si_1$  (we allow $\Si_0$ or $\Si_1$ empty).
   Let $\al\co\pi_1X\to U(n)$ be a representation.
   Fix a flat connection $B$ on $X$  with holonomy $\al$ in cylindrical
form near $\Si_0$ and $\Si_1$. The tangential  operator $A_b$ of
$D_B$ acting on $X$ decomposes as a direct sum $A_b=A_{b,0}\oplus
A_{b,1}$ since $L^2(\partial X)= L^2(\Si_0)\oplus L^2(\Si_1)$. In
particular $$\ker A_b=\ker A_{b,0}\oplus A_{b,1}\cong
H^*(\Si_0;\cc^n_\al)\oplus H^*(\Si_1;\cc^n_\al).$$
   We view $X$ as a bordism from $\Si_0$ to $\Si_1$.

   We explained in the previous section that $\ker A_b\cong H^*(\del X;
\cc^n_\al)$ is a Hermitian symplectic space. At this point we add
the hypothesis that both $\ker A_{b,0}\cong H^*(\Si_0; \cc^n_\al)$
and $\ker A_{b,1}\cong H^*(\Si_1; \cc^n_\al)$ be Hermitian
symplectic spaces. This is not automatic, but follows for example
if there exists a manifold $Y$ with boundary $\Si_0$ over which
$\al|_{\Si_0}\co\pi_1\Si_0\to U(n)$ extends.  It is in this  context
that we will usually work.

We use $X$ to define a function $L_{X,\al}$
from the set of subspaces of $H^*(\Si_0;\cc^n_\al)$ to
the set of subspaces of  $H^*(\Si_1;\cc^n_\al)$ by
\begin{equation}\label{eq2.1} L_{X,\al}(W)= P_1 \big(V_{X,\al}\cap (W\oplus
H^*(\Si_1;\cc^n_\al))\big),
\end{equation} where $P_1\co H^*(\del X;\cc^n_\al) \to
H^*(\Si_1;\cc^n_\al)$ denotes the projection onto the second factor:
$$ P_1\co H^*(\del X;\cc^n_\al)=H^*(\Si_0;\cc^n_\al)\oplus
H^*(\Si_1;\cc^n_\al) \to  H^*(\Si_1;\cc^n_\al).$$
In the following theorem, let $Y$ be a Riemannian manifold with
boundary $\Si_0$ with a product metric $g_0+ du^2 $ near the
collar.
   Write
$$Z=Y\cup_{\Si_0} X$$  and assume that $\al$ extends over $Z$. Let
$\ga_0$ be the restriction of $\ga$ to $\ker A_{b,0}$. Notice that
$\ga_0(V_{Y,\al})\oplus L_{X,\al}(V_{Y,\al})$ is a Lagrangian
subspace of $\ker A_b$.

\begin{thm}\label{lem2.1} The function of Equation \eqref{eq2.1} takes
Lagrangian subspaces to Lagrangian subspaces, i.e.\ it induces a function
\begin{equation*}L_{X,\al}\co\cL(H^*(\Si_0;\cc^n_\al))\to
\cL(H^*(\Si_1;\cc^n_\al)).
\end{equation*} This function has the properties:
\begin{enumerate}
\item If $Y$, $Z=Y\cup_{\Si_0} X$ are as above then
$$V_{Z, \al}=L_{X,\al}(V_{Y,\alpha}).$$ In short, the bordism
propagates the distinguished Lagrangian. Moreover
\begin{equation*}
\eta(D_{B, V_{Z,\al}}, Z)= \eta(D_{B,V_{Y,\al}}    ,Y) +
\eta(D_{B,\ga_0(V_{Y,\al})\oplus L_{X,\al}(V_{Y,\al})}
,X)
\end{equation*}
and hence $\rho(Z,\al,g_1)$ equals
\begin{equation*}
   \rho(Y,\al,g_0) + \rho(X,\al,g_0\amalg g_1,
\ga_0(V_{Y,\al})\oplus
L_{X,\al}(V_{Y,\al}),\ga_0(V_{Y,\tau})\oplus
L_{X,\tau}(V_{Y,\tau}) ).
\end{equation*}
where $g_i$ is a metric on $\Si_i$.
\item If $X_1$ is a bordism from $\Si_0$ to $\Si_1$ and $X_2$ is a bordism
from $\Si_1$ to $\Si_2$  and $\al\co\pi_1(X_1\cup_{\Si_1} X_2)\to
U(n)$ then $$L_{X_1\cup_{\Si_1} X_2, \al}=L_{X_2,\al}\circ L_{X_1,\al}.$$
\end{enumerate}
\end{thm}
\begin{proof}  The map of \eqref{eq2.1}  is just the map taking
$V_{X,\al}$ to its symplectic reduction with respect to the
subspace $W\oplus H^*(\Sigma_1;\cc^n_\al) \subset
H^*(\Si_0;\cc^n_\al)\oplus H^*(\Si_1;\cc^n_\al)$ (cf.\ \cite[Sec.\
6.3]{Kirk-Lesch}).   Symplectic reduction takes Lagrangians to
Lagrangians.

To prove the the first part of (1) consider a $\xi\in
V_{Y\cup_{\Si_0}X,\al}$. Then there is a $w\in
H^*(Y\cup_{\Si_0}X;\cc^n_{\al})$ with $i_{\Si_1}^*w=\xi$.
We put $\xi_0:=-i_{\Si_0}^*w$. Since certainly $w_{|X} \in
H^*(X;\cc^n_\al)$ we infer $\xi_0\oplus\xi=i_{\partial X}^*w\in
   V_{X,\al}$. Thus $\xi=P_1(\xi_0\oplus \xi)\in P_1(V_{X,\al}\cap(
V_{Y,\alpha}\oplus H^*(\Si_1;\cc^n_\al)))$.

Conversely, let $\xi\in P_1(V_{X,\al}\cap(
V_{Y,\alpha}\oplus H^*(\Si_1;\cc^n_\al)))$ be given. Then there is $\xi_0\in
V_{Y,\alpha}$ such that $\xi_0\oplus \xi\in V_{X,\al}$. Thus we may
choose $w_X\in H^*(X;\cc^n_\al)$ with $i_{\partial X}^*w_X
=\xi_0\oplus \xi$ and $w_Y\in H^*(Y;\cc^n_\alpha)$ with
$i_{\Si_0}^*w_Y=\xi_0$.

  From the Mayer--Vietoris sequence of $Y\cup_{\Si_0}X$ we obtain an
$w\in H^*(Y\cup_{\Si_0}X;\cc^n_{\al})$ with $w_{|Y}=w_Y$
and $w_{|X}=w_X$. Then $\xi=i_{\Si_1}^*w_X=i_{\Si_1}^*w\in
V_{Y\cup_{\Si_0}X,\al}$ and we reach the conclusion.

Consider now the second part of (1). We have explained in
\cite[Sec.\ 7]{Kirk-Lesch} that the gluing formula  for $\eta$--invariants
remain true if one glues (a finite union of) components of the
boundary and fixes a boundary condition at the remaining
components. The result now follows from
$V_{Z,\al}=L_{X,\al}(V_{Y,\al})$ and Theorem \ref{thm1.1}.

The proof of (2) proceeds along the same lines as the proof of the
first part of (1). Consider $W\subset H^*(\Sigma_0;\cc^n_\al)$ and
a $\xi\in L_{X_1\cup_{\Sigma_1}X_2,\al}(W)$. Then there is a
$w\in H^*(X_1\cup_{\Sigma_1}X_2;\cc^n_\al)$ with $i_{\Sigma_2}^*w=\xi$ and
$\xi_0:=-i_{\Sigma_0}^*w\in W$. Put $\xi_1:=i_{\Sigma_1}^*w$. Then
it is immediate that $\xi_0\oplus\xi_1\in V_{X_1,\al}\cap(W\oplus
H^*(\Sigma_1;\cc^n_\al))$ and $-\xi_1\oplus \xi\in V_{X_2,\al}\cap(
L_{X_1,\al}(W)\oplus H^*(\Sigma_2;\cc^n_\al))$. This proves
$\xi\in L_{X_2,\al}\circ L_{X_1,\al}(W)$.

Conversely, let $\xi\in L_{X_2,\al}\circ L_{X_1,\al}(W)$ be given. Then
there exists a $\xi_1\in H^*(\Sigma_1;\cc^n_\al)$ such that
$-\xi_1\oplus \xi\in V_{X_2,\al}\cap(L_{X_1,\al}(W)\oplus
H^*(\Sigma_2;\cc^n_\al))$
and a $\xi_0\in H^*(\Sigma_0;\cc^n_\al)$ such that
$-\xi_0\oplus \xi_1\in V_{X_1,\al}\cap(W\oplus H^*(\Sigma_1,\cc^n_\al))$.

A Mayer--Vietoris argument as in the proof of the first part of (1)
shows the existence of a $w\in H^*(X_1\cup_{\Sigma_1}X_2;\cc^n_\al)$ such
that $i^*_{\Sigma_2}w=\xi$ and $i^*_{\Sigma_0}w=-\xi_0$. This
proves $-\xi_0\oplus \xi\in V_{X_1\cup_{\Sigma_1}X_2,\al}\cap(
W\oplus H^*(\Sigma_2;\cc^n_\al))$, and hence
$L_{X_2,\al}\circ L_{X_1,\al}(W)\subset L_{X_1\cup_{\Sigma_1}X_2,\al}(W)$.
\end{proof}

   Theorem \ref{lem2.1} easily extends to the situation
$$Z=Y\cup_{\Si_0} X_1\cup_{\Si_1} \cdots \cup_{\Si_n}X_{n+1}.$$
This gives a useful strategy for computing $\rho$--invariants by
decomposing a closed manifold into a sequence of bordisms, e.g.\ by
cutting along level sets of  a Morse function.

   We use Theorem \ref{lem2.1} and the definitions    to write down a formula
   which expresses the dependence of
$\rho(Y,\al, g)$  on the metric $g$ on $\del Y$.

\begin{cor}\label{cor2.1}  Let $Y$ be a compact manifold with boundary
$\Si$. Let
$\al\co\pi_1Y\to U(n)$ be a representation. Suppose that $g_0$,  $g_1$ are
two Riemannian metrics on $\Si$. Choose a path of metrics from $g_0$ to
$g_1$ and view this path as a metric on $\Si\times [0,1]$.

Then
\[
  \begin{split} &\rho(Y,\al,g_1)-\rho(Y, \al, g_0)\\
   &=\eta(D_{B,\ga_0(V_{Y,\al})\oplus V_{Y,\al}}, \Si\times[0,1])
-\eta(D_{\Theta,\ga_0(V_{Y,\tau})\oplus
V_{Y,\tau}},\Si\times[0,1])\\
&=
\rho(\Si\times[0,1],g_0\amalg g_1,\ga_0(V_{Y,\al})\oplus
V_{Y,\al},\ga_0(V_{Y,\tau})\oplus V_{Y,\tau} ).
\end{split}\]
   \end{cor}
Here, $\Sigma$ is oriented such that a collar of the boundary takes
the form $\Sigma\times(-\ep,0]$.
\begin{proof}
Apply Theorem \ref{lem2.1} with $X=\Si\times [0,1]$ and note that
for the cylinder $X=\Si\times [0,1]$ the map $L_{X,\al}$ is the identity.
\end{proof}

In  Section \ref{Dependence on the metric} we will use Corollary
\ref{cor2.1} to
give
examples that show that $\rho(Y,\al,g)$ depends in general on the
choice of $g$, in contrast with the $\rho$--invariant for closed
manifolds.

\section{Diffeomorphism properties}

We next tie topology to the issue of the dependence of the $\rho$--invariant
on the Riemannian metric on the boundary by exploiting
the isometry invariance of spectral invariants.

Let $X$ be a manifold with boundary $\Si$.  Suppose  we are given
a diffeomorphism $f\co\Si\to \Si$ which extends to  a diffeomorphism
$F\co X\to X$. We may assume that $F$ preserves a collar of the boundary.
Then $F$ can be used to pull back flat connections and
metrics. Moreover, the diagram
\[\begin{diagram}\dgmag{900}\dgsquash
\node{H^*(X;\cc^n_\al)}\arrow{e}\arrow{s}\node{H^*(\Si;\cc^n_\al)}\arrow{s}\\
\node{H^*(X;\cc^n_{F^*(\al)})}\arrow{e}\node{H^*(\Si;\cc^n_{f^*(\al)})}
\end{diagram}
\]
shows that $f^*(V_{X,\al})=V_{X,F^*(\al)}$. Pulling back a metric
on $X$ via $F$ gives an isometry which  induces a  unitary
transformation on $L^2(\Omega^{ev}_X(E))$ taking $D_{B}$ to
$D_{F^*(B))}$ and taking $F_b^+\oplus V_{X,\al}$ to $F_b^+
\oplus V_{X,F^*(\al)}$. Hence the spectra of the two operators are
the same. Therefore
   $$\eta(D_{B,V_{X,\al}},X)=\eta(D_{F^*(B), V_{X,F^*(\al)}}, X).$$
Applying this formula first with $B$ a flat connection with holonomy $\al$ and then with $B$ the trivial connection  (and using the definition \eqref{eq1.4}) we immediately conclude the
following.
\begin{thm}\label{pii} Let $X$ be a manifold with boundary $\Si$ and
$\al\co\pi_1X\to U(n)$ a representation. Let $g$ be a Riemannian
metric on $\Si$.
   Suppose that $F\co X\to X$ is a diffeomorphism. Let $f\co\Si\to \Si$
   denote its restriction to $\Si$.
Then 
\[\rho(X,\al,g)=\rho(X,F^*(\al), f^*(g)).\tag*{\qed}\]
\end{thm}

   As an example,
recall that two diffeomorphisms $f_0,f_1\co\Si\to \Si$ are {\em
pseudo-isotopic} if there exists a diffeomorphism
$F\co\Si\times[0,1]\to \Si\times [0,1]$  which restricts to $f_0$
and $f_1$ on the boundary. In particular an isotopy is a
level-preserving pseudo-isotopy.

\begin{cor}\label{pseudo}  The
$\rho$--invariant for a manifold with boundary depends on the
Riemannian metric on the boundary only up to its pseudo-isotopy
class. Precisely,
if $f_0, f_1\co\del X\to \del X$ are  pseudo-isotopic  diffeomorphisms, then
$$
\rho(X,\al,f^*_0(g))=\rho(X,\al,f^*_1(g)).$$
\end{cor}
\begin{proof}
We may assume $f_0$ is the identity by replacing $f_1$ by $f_0^{-1}\circ f_1$ and
$g$ by $f_0^*g$. If $f_0$ is the identity then using a collar  we see that
$f_1$  extends to a diffeomorphism  $F\co X\to X$ which induces the
identity on $\pi_1X$. Thus $F^*(\al)=\al$ and  the claim follows
   from   Theorem
\ref{pii}.
\end{proof}

As another application, suppose that $F\co X\to X$ is a
diffeomorphism whose restriction to the boundary is the identity
(or pseudo-isotopic to the identity). Then for any metric $g$ on
$\Si$ $$\rho(X,\al,g)=\rho(X,F^*(\al),g).$$
These facts can be summarized as follows.
   Let $\mathscr{M}_\Si$
denote the space of Riemannian metrics on $\Si$, and  $\cF_X$ the
space of flat connections on $E\to X$.  If $\del X=\Si$, let
$\cD_X$ denote the diffeomorphism group of $X$ and let $\cD_{X}^0$
denote subgroup of those diffeomorphisms which induce the identity
on $\pi_1(X)$. This group acts on $\cM_\Si$.

The assignment $$(g,B)\mapsto \rho(X,\al,g)$$ defines  a function
$$\rho(X)\co\cF_X\times\mathscr{M}_\Si\to \rr.$$
   Theorems \ref{thm1.1}   and  \ref{pii} say that $\rho(X)$
    descends to a function on the quotient
\begin{equation}\label{rhofunct}
\rho(X)\co\chi(\pi_1X,U(n))\times(\mathscr{M}_\Si/\mathscr{D}^0_X)\to
\rr,\end{equation}
   where
\begin{equation}\label{defofchi}
\chi(\pi_1X,U(n))=\Hom(\pi_1X,U(n))/\text{conj.}=\cF_X/\cG_X
\end{equation}
with $\cG_X$ the group of gauge transformations of $E\to X$.  The
quotient $\cD_X/\cD^0_X$ acts diagonally on
$\chi(\pi_1X,U(n))\times(\mathscr{M}_\Si/\mathscr{D}^0_X)$ and the
function of \eqref{rhofunct} is invariant under this action.

\section{Dependence on the metric}\label{Dependence on the metric}

In this section we prove the following theorem which shows that the
$\rho$-
and $m$--invariants  depend on the choice of metric on the boundary.

\eject

\begin{thm} \label{dependence}\indent\par
\begin{enumerate}
\item There exists a 3-manifold $Y$ with non-empty boundary,  Riemannian
metrics $g_0$, $g_1$ on $\del Y$, and a representation $\al\co\pi_1Y\to
U(2)$ so that
$$\rho(Y,\al,g_0)\ne \rho(Y,\al ,g_1).$$ Examples exist with vanishing
kernel of the tangential operator, i.e.\ 
\[\ker A_b\cong H^*(\del Y;\cc^2_\al)=0.\]
\item There exist metrics $g_0$ and $g_1$ on the torus $T$ and 3-manifolds
$X$ and $Y$ with boundary $T$  such that setting
$V_X=  \text{image}\big( H^*(X;\cc)\to H^*(T;\cc)\big)$ and
$V_Y= \text{image}\big(  H^*(Y;\cc)\to H^*(T;\cc)\big)$  (with
$\theta$ the trivial connection),
$$m(V_X,V_Y)_{(\theta,g_0)}\ne m(V_X,V_Y)_{(\theta,g_1)}.$$
\end{enumerate}
\end{thm}

\begin{sloppypar}
In the first statement of Theorem \ref{dependence} the point of taking an
example with
$H^*(\del Y;\cc^n_\al)=0$ is to emphasize that the metric dependence of
the $\rho$--invariant is much more subtle than just being a consequence of
the dependence of $m(V,W)$ on the metric.
\end{sloppypar}

To understand the significance of the second statement, observe that the
choice of Riemannian metric on $\del X=\Sigma$ enters into the definition
of $m(V,W)$ only through the restriction of the induced $L^2$ metric on
$L^2(E|_\Sigma)$ to the harmonic forms
$\ker A_b\cong H^*(\Sigma;\cc^n)$.  There are clearly many Riemannian
metrics on $\Sigma$ which restrict to the same metric on the space of
harmonic forms.
It is perhaps at least intuitively clear that the invariant $m(V,W)$ of
pairs of Lagrangian subspaces in a   Hermitian
symplectic  space can vary   as the inner product varies. But
our argument shows    more: the metrics we use are   restrictions of
$L^2$ metrics to the harmonic forms (i.e.\ the kernel of the tangential
operator) and the Lagrangian subpaces we consider are of the form
image$\big(H^*(X;\cc^n_\al))\to H^*(\del X;\cc^n_\al)\big)$.    Notice that
these Lagrangians $V_X$ are always graded direct sums; i.e.\
$V_X=\oplus_i V_X^i$ with $V_X^i=$ image$\big(H^i(X;\cc^n_\al)\to H^i(\del
X;\cc^n_\al)\big)$.

As an illuminating  non-example the reader might consider the case when
$\del X$ is a $2k$-sphere, and $\al$ is trivial. Then
$\ker A_b=H^*(S^{2k})=H^0(S^{2k})\oplus H^{2k}(S^{2k})$. Certainly one
can find   families of Riemannian metrics on
$S^{2k}$ so that the induced metric on the harmonic forms $H^*(S^{2k})$
varies (e.g.\ by scaling the metric) and from that it is not hard to produce
a pair of Lagrangian subspaces
$V,W\subset H^*(S^{2k})$ for which
$m(V,W)_{(\theta,g)}$ varies with $g$. But, if $\al$ is the trivial
representation on $\pi_1X$, then (for any such $X$)  the subspace
$V_X=\text{image }H^*(X;\cc^n_\al)\to H^*(\del X;\cc^n_\al)$ is just
$H^0(S^{2k})$.  Therefore,  given a similar
$Y$,
$V_X=V_Y$ and so $\phi(V_X)\phi(V_Y)^*=Id$. This implies that the {\it
geometric} invariant $m(V_X,V_Y)_{(\theta,g)}$ is independent of $g$ in this
situation.

\subsection{$\rho(X,\al,g)$ depends on $g$}

We begin with  the proof of the first part of Theorem \ref{dependence} by
providing an explicit example which shows that
$\rho(Y,\al,g)$     depends in general on the choice of Riemannian metric
$g$ on the boundary $\del Y$.   We will show that there exists a
3-manifold $Y$ with boundary a torus, a non-abelian representation
$\al\co\pi_1Y\to SU(2)$ with $H^*(\del Y;\cc_\al^2)=0$, and   Riemannian
metrics $g_0$ and $g_1$ on the torus so that
$\rho(Y,\al,g_1)-\rho(Y,\al,g_0)$ is non trivial. The metrics $g_0$ and
$g_1$ can be taken to be flat.

The manifold $Y$ we take is the complement of the right-handed trefoil
knot in
$S^3$.  The analysis of the space of $SU(2)$ representations of the
fundamental groups of knot complements has a long history in the
literature, starting with the beautiful article
\cite{klassen-thesis}.  Details and proofs of most of the facts we use
here can be found in
\cite{BHKK}.

The fundamental group of
$Y$ is
$$\pi_1 Y=\langle x,y\ | \ x^2=y^3\rangle.$$ The boundary of $Y$ is a
torus, and the  meridian $\mu$ and longitude $\la$ of $Y$ generate
$\pi_1(\del Y)=\zz^2.$ They are   given  in this presentation of
$\pi_1Y$ by $$\mu=xy^{-1}\text{ and }\la=x^2(xy^{-1})^{-6}.$$
The space of conjugacy classes of non-abelian $SU(2)$ representations of
$\pi_1(Y)$ is an open arc. Moreover, given any pair $(\phi,\psi)$ in the
open line segment in $\rr^2$
\begin{equation}\label{open arc}\{(t, -6t+\tfrac{1}{2})\ | \
\tfrac{1}{12}<t<\tfrac{5}{12}\} \end{equation} there exists  a unique
conjugacy class of non-abelian $SU(2)$ representations  of
$\pi_1(Y)$ which satisfies
\begin{equation}\label{eq3.3}\mu\mapsto
    \begin{pmatrix} e^{2\pi i\phi }&0\\0&e^{-2\pi i\phi }\end{pmatrix},\
\la\mapsto \begin{pmatrix} e^{2\pi i\psi }&0\\0&e^{-2\pi i\psi
}\end{pmatrix}.
\end{equation}
Therefore, letting
\begin{equation}\label{eq3.6}(\phi_1,\psi_1)=(\tfrac{1}{5},-\tfrac{7}{10})
\text{ and }(\phi_2,\psi_2)=(\tfrac{2}{5},-\tfrac{19}{10}),\end{equation}
we obtain two non-abelian representations $\al_1,\al_2\co\pi_1Y\to SU(2)$
in the open arc of \eqref{open arc}.

Fix an identification of the boundary of $Y$ with the 2-torus
$T=\rr^2/\zz^2$ such that $\mu$ corresponds to the $x$-axis and $\la$ to
the $y$-axis.  Give $T$ the induced flat metric $g_0$.

   Consider the matrix
\begin{equation}\label{eq3.14}f=\begin{pmatrix} 1&3\\2&7\end{pmatrix}\in
SL(2,\zz).\end{equation} Then $f$ acts by right multiplication on $\rr^2$
preserving the standard lattice
$\zz^2$, and hence induces a diffeomorphism
$f\co T\to T$.

The first part of Theorem \ref{dependence} follows from the next theorem.

\begin{thm}\label{thm3.1} The difference
$\rho(Y,\al_1,g_0)-\rho(Y,\al_1,f^*(g_0))$  does not equal the difference
$\rho(Y,\al_2,g_0)-\rho(Y,\al_2, f^*(g_0))$. Hence the $\rho$--invariant
for manifolds with boundary depends in general on the choice of
Riemannian metric on the boundary, and moreover
   $\rho(Y,\al,g)-\rho(Y,\al,f^*(g))$ is  not a function of $f\co\del Y\to
\del Y$
   alone.
\end{thm}

\begin{proof}

Let $\hal_i\co\pi_1T\to SU(2)$ denote the restrictions of $\al_i$ to
$\pi_1(\del Y)$. Let  $g_1$ denote the pulled back metric $g_1=f^*(g_0)$.

Fix a Riemannian metric $\tilde{g}$ on $Y$ in product form near the
boundary so that the restriction of $\tilde{g}$ to the boundary equals
$g_0$.

Choose a smooth path
$g_t$ of Riemannian metrics on $T$   from $g_0$ to $g_1$ which is
stationary for $t\in [0,\ep]$ and $[1-\ep,1]$. Then $g_t$     determines
the    metric $g_t + dt^2$  on $T\times [0,1]$.

Let $B_1$ be a flat  connection on $Y$  with holonomy $\al_1$ and in
cylindrical form $B=\pi^*(b_1)$ near $\del Y$. Let $D_{B_1}$ denote the
corresponding odd signature operator on $Y$. Then $D_{B_1}$ has an
obvious extension to
$Y\cup(T\times [0,1])$ by defining it to be the pullback of $b_1$ via the
projection $T\times [0,1]\to T$.  Similarly choose a flat connection
$B_2$ with holonomy $\al_2$ and extend it to $ T\times [0,1]$.

\begin{lemma}\label{lem3.1} $H^*(T;\cc^2_{\hal_i})=0$.
\end{lemma}

\begin{proof} Applying  the Fox calculus  to the presentation
$$\pi_1(\del Y)=\langle \mu, \la \ | \ \mu\la\mu^{-1}\la^{-1}\rangle$$ we
conclude that  $H^*(T;\cc^2_{\hal_i}) $ is the cohomology of the complex
$$0\to \cc^2\xrightarrow{\del_0}\cc^2\oplus
\cc^2\xrightarrow{\del_1}\cc^2\to 0,$$ where
$$\del_0=\begin{pmatrix}\hal_i(\mu)-I & \hal_i(\la)-I\end{pmatrix}
\text{ and }
\del_1=\begin{pmatrix}I-\hal_i(\la) \\ \hal_i(\mu)-I\end{pmatrix}.$$
A simple computation using \eqref{eq3.3}  and \eqref{eq3.6} shows 
that the cohomology of this complex vanishes.
\end{proof}

Continuing with the proof of Theorem \ref{thm3.1}, It follows from
Corollary \ref{cor2.1} that
\begin{equation}\label{eq3.7}
\begin{split}
\big( \rho(Y, \al_1, g_1) &- \rho(Y,\al_1, g_0)\big) -\big( \rho(Y, \al_2,
g_1) - \rho(Y,\al_2, g_0)\big) \\ 
&= \eta(D_{B_1},\Si\times[0,1]) -\eta(D_{B_2},\Si\times[0,1])
\end{split}\end{equation} (Note: Lemma \ref{lem3.1} implies that $\ker
A_{b_i}=0$ so that there are no Lagrangian subspaces to specify in the
$\eta$--invariants in
   \eqref{eq3.7}.)

We will show that the right side of \eqref{eq3.7} is not an integer, from
which Theorem \ref{thm3.1} follows.

Let $M_f$ denote the mapping torus of $f$:
$$M_f=T\times[0,1]/(t,0)\sim (f(t),1).$$ The metric on $T\times [0,1] $
descends to a metric on
$M_f$  since the gluing map
$f\co(T,g_1)\to (T,g_0)$ is an isometry.

Recalling that  $\mu$ and $\la$ in $\pi_1(T)$ denote the two generators,
$$\pi_1(M_f)=\langle \mu,\la,\tau \ | \ [\mu,\la]=1,
\tau\mu \tau^{-1}= \mu\la^2, \tau\la
\tau^{-1}=\mu^3\la^7\rangle.$$
 It follows that given a pair of real numbers
$\phi,\psi$, the assignment
\begin{equation}\label{eq3.2}\mu\mapsto \begin{pmatrix} e^{2\pi
i\phi}&0\\0&e^{-2\pi i\phi}\end{pmatrix},\
\la\mapsto \begin{pmatrix} e^{2\pi i\psi}&0\\0&e^{-2\pi
i\psi}\end{pmatrix}
\text{ and }
\tau\mapsto \begin{pmatrix} 0&1\\-1&0\end{pmatrix}\end{equation}
determines a representation
$ \pi_1(M_f)\to SU(2) $ if and only if
$e^{-2\pi i\phi}=e^{2\pi i(1\phi+2\psi)}$ and
$e^{-2\pi i\psi}=e^{2\pi i(3\phi+7\psi)}$, i.e.~if and only if
\begin{equation}\label{eq3.1}\begin{pmatrix}\phi&
\psi\end{pmatrix}(f_*+\text{Id})\equiv 0\pmod{\zz}.\end{equation}
Equation \eqref{eq3.1} holds for $(\phi_1,\psi_1)$ and $(\phi_2,\psi_2)$
as in
\eqref{eq3.6}.

    Thus taking $(\phi_1,\psi_1)$ and $(\phi_2,\psi_2)$ in \eqref{eq3.2} we
obtain two representations
$$\beta_j\co\pi_1(M_f)\to SU(2),\  j=1,2,$$ with the property that their
restrictions to the fiber $T\times \{0\}$  equal
$\hal_j$.

View $M_f$ as the union of two cylinders $T\times [0,1] \cup T\times
[0,1]$ using the gluing map $Id \cup f$.  Give $M_f$ the product metric
$g_0 + dt^2$ on the first piece and $g_t + dt^2$ on the second. Equation
   \eqref{eq1.49}  shows that
\begin{equation}\label{eq3.8} \eta(D_{B_j}, M_f) = \eta(D_{B_j},
(T\times[0,1], g_0 + dt^2)) +
\eta(D_{B_j}, (T\times[0,1], g_t + dt^2)).\end{equation}
It follows from Corollary \ref{cor2.1} that $$ \eta(D_{B_1},
(T\times[0,1], g_0 + dt^2))- \eta(D_{B_2}, (T\times[0,1], g_0 +
dt^2))=0$$ (alternatively
   Lemma 7.1 of \cite{Kirk-Lesch} shows directly that
$ \eta(D_{B_j}, (T\times[0,1], g_0 + dt^2))=0$).

Thus combining
\eqref{eq3.8} for $j=1,2$ with \eqref{eq3.7} we obtain
\begin{equation}\label{eq3.9}
\begin{split}
\big( \rho(Y, \al_1, g_1) &- \rho(Y,\al_1, g_0)\big) -\big( \rho(Y, \al_2,
g_1) - \rho(Y,\al_2, g_0)\big)\\ &= \eta(D_{B_1}, M_f) - \eta(D_{B_2},
M_f)\\
&=\rho(M_f, \be_1)-\rho(M_f,\be_2).\end{split}\end{equation} The
last equality follows from the definitions of $\eta$ and $\rho$ for a
closed manifold.

We have thus reduced the problem to showing that the difference of the
$\rho$--invariants for $\be_1$ and $\be_2$ on the {\em closed} manifold
$M_f$ is not an integer. On a closed manifold, the $\rho$--invariants and
the Chern-Simons invariants are related by the Atiyah-Patodi-Singer
theorem \cite{APS-I,APS-II}; the formula is
   (see \cite[Sect.\ 5.3-5.5]{BHKK}):
\[\begin{split}
\SF(D_{B_t}) = 2(&\text{cs}(B_1) - \text{cs}(B_2))\\
    &+\tfrac{1}{2}(\rho(M_f, \be_1)-\rho(M_f,\be_2) 
         -\dim\ker D_{B_1}+\dim \ker D_{B_2}).
  \end{split}
\]
Here $\SF(D_{B_t})$ denotes the spectral flow (an integer) of
the family of self-adjoint elliptic operators $D_{B_t}$ where $B_t$ is
any family of connections from $B_1$ to $B_2$. This implies
\begin{equation}\label{eq3.10}
\rho(M_f, \be_1)-\rho(M_f,\be_2)\equiv
4(\text{cs}(B_2)-\text{cs}(B_1))\pmod{\zz}.\end{equation}
Theorem 5.6 of \cite{Kirk-Klassen:chern-simons} calculates the
Chern-Simons invariant mod $\zz$ of  flat connections on $M_f$ in terms
of the vector
$(\phi,\psi)$ and the matrix $f$: if
$(m,n)=(\phi,\psi)(\text{I}+f^{-1})$, then the Chern-Simons invariant of
the flat connection with holonomy representation determined by
$(\phi,\psi)$ equals
$\phi n-\psi m$ mod
$\zz$.

    Since $(\phi_1,\psi_1)(I+f^{-1})=(3,-2)$ and
$(\phi_2,\psi_2)(I+f^{-1})=(7,-5)$ this gives
\begin{equation} cs(B_1)\equiv
\frac{7}{10}\pmod{\zz}\text{ and }cs(B_2)\equiv
\frac{3}{10}\pmod{\zz}.\end{equation} Hence
\begin{equation}\label{eq3.11}
4(\text{cs}(B_1)-\text{cs}(B_2))\equiv\tfrac{28}{10}
-\tfrac{12}{10}\equiv\tfrac{3}{5}\pmod{\zz}.\end{equation}
Combining  \eqref{eq3.9}, \eqref{eq3.10}, and \eqref{eq3.11}
   we see that
$$\big( \rho(Y, \al_1, g_1) - \rho(Y,\al_1, g_0)\big) -\big( \rho(Y,
\al_2, g_1) - \rho(Y,\al_2, g_0)\big)\ne 0,
$$
proving Theorem \ref{thm3.1} and hence the first assertion of Theorem
\ref{dependence}.\end{proof}

An interesting problem suggested by Corollary \ref{cor2.1} and  Theorem  \ref{thm3.1} is to
find a description of the function $\chi(\Si, U(n))\times
\cD_\Si/\cD^0_\Si\to\rr$ which takes $(\al, f)$ to
$\rho(\Si\times[0,1],\al, g \cup f^*(g))$ (i.e.~the $\rho$--invariant of
the cylinder with a fixed metric $g$ at $\Si\times\{0\}$ and
$f^*(g)$ at $\Si\times\{1\}$). Theorem \ref{thm3.1} implies that
this map is non-trivial, and depends on $\al$.

\subsection{$m(V_{X,\al},V_{Y,\al})_{(b,g)}$ depends on
$g$}
We next prove the second assertion of Theorem \ref{dependence}.

Consider the 2-torus $T^2=S^1\times S^1$ with its standard oriented basis
of 1-forms $\{dx, dy\}$. We consider these forms as  sections of the
trivial 1-dimensional complex bundle over $T^2$ endowed with the trivial
connection.

For each $t>0$, give $T^2$ the Riemannian metric  for which
$\{dx, t\ dy\}$ is an orthonormal basis at each point. Letting $\hat{*}_t$
denote the corresponding Hodge $*$-operator we have
$$\hat{*}_t dx= t\ dy, \ \ \hat{*}_t dy=-\tfrac{1}{t}\ dx,\ \ \hat{*}_t 1
= t\ dx\wedge dy,
\ \
\text{ and }\ \hat{*}_t (dx\wedge dy)= \tfrac{1}{t}.$$
    Hence
$$\ga_t(dx)= t\ dy, \ \ \ga_t(dy)=-\tfrac{1}{t}\ dx,\ \ \ga_t(1) = t\
dx\wedge dy,
\ \
\text{ and }\ \ga_t(dx\wedge dy)= -\tfrac{1}{t}.$$
    This defines the de Rham   operator  $A_t=\pm(\hat{*}_td+d\hat{*}_t)$ as
above.

The harmonic forms $\cH_t^*=\ker A_t $ with respect to this metric are
independent of
$t$ as one can readily compute: the harmonic 0-forms $\cH^0$ are the
constant functions,  the harmonic 1-forms  $\cH^1$ are $a\ dx + b \ dy$
with $a,b$ constant, and the harmonic 2-forms $\cH^2$ are $a\ dx\wedge
dy$ with $a$ constant.

We can compute the $L^2$ inner product $\langle\ , \ \rangle_t$
restricted to the hamonic forms:
$$\langle 1,1\rangle_t=\int_T 1 \wedge \hat{*}_t 1=\int_T t\ dx\wedge dy
= 4
\pi^2 t.$$ Similarly
\[\begin{split}
&\langle dx,dx\rangle_t= 4\pi^2 t, \ \
\langle dx,dy\rangle_t= 0, \ \
\langle dy,dy\rangle_t= 4\pi^2 /t, \\
& \text { and }
\ \langle dx\wedge dy,dx\wedge dy\rangle_t= 4\pi^2/t.
  \end{split}
\]
The Hermitian symplectic space of harmonic forms $(\cH^*,\ga_t,\langle \
, \ \rangle_t)$ is a direct sum
$ (\cH^0\oplus \cH^2)\oplus\cH^1$ of two Hermitian symplectic spaces.
Thus there is a corresponding splitting of the  $\pm i$ eigenspaces of
$\ga_t$. One checks that the $\pm i$ eigenspaces of $\ga_t$ acting on
$\cH^*$ are (with the obvious notation)
$$ (\cH^0\oplus \cH^2)^\pm = \text{span}\{1 \mp i t \ dx\wedge dy\}
$$
and $$ (\cH^1)^\pm = \text{span}\{dx \mp i t \ dy\}.
$$
Suppose that $X$ is a compact 3-manifold with $\del X=T$. Then (taking
$\cc$ coefficients)
\[\begin{split} V_X&=\text{image}\big( H^*(X)\to H^*(T)\big)\\
& =\text{image}\big(H^0(X)\to H^0(T)\big)\oplus \text{image}
(H^1(X)\to H^1(T)\big)\oplus\\
&\quad\oplus  \text{image}\big( H^2(X)\to H^2(T)\big)\\
& =
H^0(T)\oplus
   \text{image}\big( H^1(X)\to H^1(T)\big)\oplus 0.\end{split}
\] Write $V_X^1$ for $\text{image}\big( H^1(X)\to H^1(T)\big)$.  Similarly if
$Y$ is another manifold with $\del Y=T$ we have $V_Y=H^0(T)\oplus V_Y^1$ with
$V_Y^1=\text{image}\big(H^1(Y)\to H^1(T)\big)$.

   Notice that since the coefficients are obtained by tensoring the integer
cohomology with $\cc$, there exist integers $a,b,A,B$ so that
$V_X^1=\text{span}\{a\ dx + b \ dy\}$ and
$V_Y^1=\text{span}\{A\ dx + B \ dy\}$.  Moreover, given any pair of (not
both zero) integers $(a,b)$ one can find a 3-manifold $X$ with
$V_X^1=\text{span}\{a\ dx + b
\ dy\}$.

For example,  take $X=S^1\times D^2$. By clearing denominators we may
assume that $a$ and $b$ are relatively prime. Suppose that $p,q$ are
integers satisfying
$ap-bq=1$. Then there is a diffeomorphism $\del X= S^1\times S^1$ to
$T^2$ covered by the linear map $\rr^2\to \rr^2$ with matrix
\[\begin{pmatrix} a & q\\b&p\end{pmatrix}.\] The closed 1-form $dx$ on
$\del X$ is identified with $a\ dx + b\ dy$ on $T^2$. Since $dx$ extends
to $X$, this gives an example with $V_X^1=\text{span}\{a \ dx + b \
dy\}$.

In terms of the $\pm i$ eigenspace decomposition of $\cH^*$ one can
easily check that
$$\phi(V_X)=\phi^{0,2}(H^0(T))\oplus \phi^1(V_X^1),$$ where
$$\phi^{0,2}(H^0(T))(1-it\ dx\wedge dy)=1+i t \ dx\wedge dy$$ and
$$\phi^1(V_X^1)(dx - i t \ dy) = \frac{ita+b}{ita -b}(dx + it \ dy).$$
(See Equation \eqref{eqforphi} for the definition of the unitary map $\phi(V)$ associated to a Lagrangian subspace $V$.)
These equations imply that
\[\phi(V_X)\phi(V_Y)^*=\begin{pmatrix} 1&0\\0& (\frac{ita+b}{ita
-b})(\frac{itA-B}{itA +B})\end{pmatrix}.\] Therefore (see \eqref{defofm})
\[m(V_X,V_Y)_{(\tau,g_t)}=-\tfrac{1}{\pi i}\big(\pi i + \log(-(\frac{ita+b}{ita
- b})(\frac{itA-B}{itA
+B}))\big) + \dim (V_X\cap V_Y).\]
For example, taking $B=0$ this reduces to
\[m(V_X,V_Y)_{(\theta,g_t)}=-1 + \dim(V_X\cap V_Y) -\tfrac{1}{\pi i}
\log(\frac{b+ita}{b- ita }).\]
But $\frac{b+ita}{b-ita }=\frac{(b+ita)^2}{|b+ita|^2}$ and so
$\log(\frac{b+ita}{b-ita })$ is equal to the argument of $(b+ita)^2$
which varies non-trivially as $t$ varies provided both $a$ and $b$ are
non-zero.

Thus we have given an example of a family of Riemannian metrics $g_t$ on
the torus $T$ and  shown how to find 3-manifolds $X$ and $Y$ so that
(with respect to the trivial $U(1)$ representation $\tau$)
$m(V_X,V_Y)_{(\theta,g_t)}$ varies non-trivially as $t$ is varied.  This proves
the second part of Theorem \ref{dependence}. \qed

\section{An extension of the Farber--Levine--Weinberger\newline
theorem to manifolds with boundary}\label{flwtheorem}

Suppose that $F\co M\to M'$ is an orientation preserving  homotopy equivalence of
smooth compact manifolds. Then $F$ induces an isomorphism of fundamental
groups, and hence a homeomorphism (in fact a real-analytic isomorphism)
$$\Hom(\pi_1(M'),U(n))\xrightarrow{F^*}\Hom(\pi_1(M),U(n)).$$
Taking the quotient by the action of conjugation eliminates the
dependence on base points, and   one obtains an
   identification (see Equation \eqref{defofchi})
$$\chi(\pi_1(M'),U(n))=\chi(\pi_1(M),U(n)).$$
If $M$ and $M'$ are closed, then taking $\rho$--invariants defines
functions (write $\pi=\pi_1M$ for convenience)
$$\rho(M)\co\chi(\pi,U(n))\to \rr\text{ and }
\rho(M')\co\chi(\pi,U(n))\to
\rr.$$
In \cite{Farber-Levine} M. Farber, J. Levine, and S. Weinberger proved
the following remarkable theorem.
\begin{thm}[Farber-Levine, Weinberger] \label{FLW} The difference
$$\rho(M)-\rho(M')\co\chi(\pi,U(n))\to \rr$$
factors through the set of path components of $ \chi(\pi,U(n)) $
and takes values in the rationals. Briefly, there is a commutative
diagram
\[\begin{diagram} \dgmag{1200}\dgsquash
\node{\chi(\pi,U(n))}\arrow{e,t}{\rho(M)-\rho(M')}\arrow{s}\node{\rr}\\
\node{\pi_0(\chi(\pi,U(n)))}\arrow{e}\node{\qq}\arrow{n}
\end{diagram}\]
Moreover the difference $\rho(M)-\rho(M')$ vanishes on the path
component containing the trivial representation.
\qed
\end{thm}

Their proof has 3 ingredients. First Farber and Levine show that the difference
$\rho(M)-\rho(M')$ modulo $\zz$  factors through the set of path
components using the Atiyah-Patodi-Singer theorem and a computation
of the index density. We will generalize this fact using the Dai-Freed theorem
in Theorem \ref{cor16}.

Next they show that the ``$\zz$ part'',
i.e.\ the spectral flow of the odd signature operator along a path of
flat connections on a closed manifold, is a homotopy invariant which
can be derived from a certain linking form. (A slightly different
argument for this part was given in \cite{KK3}.)

Finally in an appendix Weinberger
uses algebraic techniques to show that the difference is rational.

   In this and the following section we will extend these results
to manifolds with boundary, with respect to homotopy equivalences
which restrict to diffeomorphisms on the boundary. (One cannot hope
to prove a generalization for homotopy equivalences which do not
behave nicely on the boundary; see Theorem \ref{squaregranny}.)

Suppose that $F\co X\to X'$ is a smooth map between compact manifolds
which restricts to a diffeomorphism $f=F|_{\del X}\co\del X\cong \del X'$
on the boundary.   Pulling back representations of $\pi_1(X')$ and
Riemannian metrics on $\del X'$ induces a function (an analytic
isomorphism if $F$ induces an isomorphism on fundamental groups):
$$\chi(\pi_1X',U(n))\times \mathscr{M}_{\del X'}\to
\chi(\pi_1X,U(n))\times \mathscr{M}_{\del X}.$$
In particular if $F$ is a homotopy equivalence we consider
$\rho(X)$ and $\rho(X')$ as functions on the same space via this
identification.  Write $\pi$ for $\pi_1X$.

\begin{thm}\label{FLWKL} Let $F\co X\to X'$ be a homotopy equivalence of compact
manifolds which restricts to a diffeomorphism on the boundary.  Then
the difference
   $$\rho(X)-\rho(X')\co\chi(\pi,U(n))\times \mathscr{M}_{\del X}\to \rr$$
factors through  $\pi_0(\chi(\pi, U(n)))\times (\mathscr{M}_{\del X}/
\mathscr{D}^0_{\del X})$
(where $\mathscr{D}^0_{\del X}$ denotes the group of diffeomorphisms
of $\del X$ pseudo-isotopic to the identity)  and takes values in
the rational numbers.

In other words there is a commutative diagram
\[\begin{diagram} \dgmag{1000}\dgsquash
\node{\chi(\pi,U(n)) \times \mathscr{M}_{\del X}}\arrow{e,t}{\rho(X)-
\rho(X')}\arrow{s}\node{\rr}\\
\node{\pi_0(\chi(\pi,U(n))) \times (\mathscr{M}_{\del X}/
\mathscr{D}^0_{\del X})}\arrow{e}\node{\qq}\arrow{n}
\end{diagram}\]
Moreover the difference $\rho(X)-\rho(X')$ vanishes on the path component of
the  trivial representation.

\end{thm}
\begin{proof}
We may asume, by homotoping $F$ slightly, that $F$ restricts to a
diffeomorphim  of  collar neighborhoods  of the boundary. Identify
$\Sigma=\del X$ with $\del X'$ via $f$, and fix a metric $g$ on
$\Sigma$.

Let $\al$  be a $U(n)$ representation of $\pi_1(X)$.
Consider the  $(2k+2)$-manifold $W=X\times [0,1]$. Then
$\al$    clearly extends  to $\pi_1W$. By smoothing the
corners of
$W$ we obtain a smooth manifold with boundary $\del W= (-X) \cup_\Sigma
X$. Since $W$ is   a product, $\Sign(W)$
and $\Sign_{\al}(W)$ both vanish.  The Atiyah-Patodi-Singer
theorem  then implies that
   $\rho(\del W,\al)=0$.

(Alternatively, there is a direct
spectral argument which shows the vanishing
of $\rho(\partial W,\al)$: the reflection which interchanges the two
copies of $X$ in $\partial W=(-X)\cup_\Sigma X$ is orientation reversing, hence
it anticommutes with the odd signature operator $D_B$. Thus the
spectrum of $D_B$ is symmetric and so its $\eta$--invariant vanishes.)

The homotopy equivalence $F\co X\to X'$ induces a homotopy equivalence
of closed manifolds
$$\operatorname{Id}\cup F\co \del W=(-X)\cup_\Sigma X\to (-X)\cup_\Sigma X'.$$
The Farber--Levine--Weinberger theorem then implies that
\begin{equation}\label{eq7.3}
\rho((-X)\cup_\Sigma X',\al)=\rho((-X)\cup_\Sigma X',\al)-\rho(\del
W,\al)=r\in
\qq\end{equation}
for some rational  number $r$ which depends only on the path
component of $\al$ in $\chi(\pi_1(W), U(n))=\chi(\pi,U(n))$.

Using Theorem \ref{thm1.1}, part 4 we conclude
that
\begin{equation}
\label{eq7.1}
   \rho((-X)\cup_\Sigma X',\al)=
      \rho(X', \al, g) + \rho(-X, \al ,g,
\ga(V_{X',\al}),\ga(V_{X',\tau}))
   \end{equation}
and
\begin{equation}
\label{eq7.2}
   \rho((-X)\cup_\Sigma X,\al)=
      \rho(X, \al, g) + \rho(-X, \al ,g,
\ga(V_{X,\al}),\ga(V_{X,\tau})).
   \end{equation}
The commutative diagram (with any coefficients)
\[\begin{diagram}\dgmag{800}
\node{H^*(X')}\arrow[2]{e}\arrow{se,b}{F^*}\node[2]{H^*(\Sigma)}\\
\node[2]{H^*(X)}\arrow{ne}
\end{diagram}\]
shows that $V_{X',\al}=V_{X,\al}$ and  $V_{X',\tau}=V_{X,\tau}$.
Therefore,
\begin{equation}\label{eq7.4}\rho(-X, \al ,g,
\ga(V_{X,\al}),\ga(V_{X,\tau}))=\rho(-X, \al ,g,
\ga(V_{X',\al}),\ga(V_{X',\tau})).\end{equation}
Taking the difference of \eqref{eq7.1} and \eqref{eq7.2} and using
\eqref{eq7.3} and \eqref{eq7.4} we conclude that
$$ \rho(X',\al,g) -\rho(X, \al, g)=\rho((-X)\cup_\Sigma X',\al)=r\in \qq$$
for an $r$ that depends only on the path component of $\al$ in
$\chi(\pi,U(n))$.  Notice that if $\al$ is trivial, then $\rho(X,\al,g)=0$.

The fact that $\rho(X,\al,g)$ depends only on the pseudo-isotopy class
of $g$ follows from Corollary \ref{pseudo}.
\end{proof}

The reader should keep in mind that the $\rho$--invariants in the context of
Theorems \ref{FLW} and \ref{FLWKL} are not continuous in $\al$. This is because
eigenvalues of the odd signature operator can become zero, or change sign, as
$\al$ is varies.   Thus what is being asserted in these theorems is that the
discontinuities of $\rho$ are homotopy invariants, provided the homotopy equivalence
restricts to a diffeomorphism on the boundary.

We formalize and extend this remark in Theorem  \ref{sfthm} below  which shows
that the spectral flow of the odd signature operator coupled to a path of flat
connections on a manifold with boundary is a homotopy invariant. For a closed
manifold this is the main result of
\cite{Farber-Levine}, and the principal ingredient in the proof of Theorem
\ref{FLW}.
Partial results for manifolds with boundary were obtained in a series of
articles by E. Klassen and the first author, including \cite{KK1, KK2, KK4}  as
well as in the articles \cite{KKR, BHKK} which also contain applications of
these ideas to   calculations of Floer homology, $SU(3)$ Casson invariants,
and TQFT.

Consider, then, a path $B_t, t\in [0,1]$ of flat $U(n)$ connections, in
cylindrical form near the boundary, on a compact smooth manifold $X$ with
boundary $\del X$.  These give a path $D_{B_t}$ of odd signature operators in
the form $D_{B_t}=\ga(\tfrac{\del}{\del x} + A_{b_t})$ on a collar neighborhood
of $\del X$.

To obtain a path of {\em self-adjoint} operators on $X$ whose kernels have
a topological meaning we assume that we are given a continuous path of APS
boundary conditions.   Precisely, we assume  that we are given a
continuous path
$P(t)$ of Lagrangian subspaces of $L^2(E|_{\del X})$ so that for each $t$
there exists a (finite-dimensional) subspace $W(t)\subset \ker A_{b_t}$ for
which
$$P(t) = F_{b_t}^+(t)\oplus W(t).$$
The path of operators $D_{B_t,P(t)}$ (recall this means $D_{B_t}$ with the
boundary condition given by $P(t)u=0$)
   is a path of self-adjoint discrete operators (see  Section
\ref{Section2}, \cite{Kirk-Lesch}, and in particular
\cite[Sec.\ 3]{BooLesPhiUFS}) and hence  has a spectral flow
$\SF(X,D_{B_t,P(t)})_{t\in [0,1]}\in \zz$.
We use the $(-\epsilon,-\epsilon)$ convention for spectral flow; this implies
that the spectral flow is additive with respect to composition of paths.

To ensure that the spectral flow of the resulting path of
self-adjoint operators
$D_{B_t,P(t)}$  is a topological invariant we furthermore assume that
  $$P(0)=F^+_{b_0}\oplus V_{X,\al_0}  \text{ and } P(1)=F^+_{b_1}\oplus
V_{X,\al_1}.$$
As before $F_{b_t}^+(t)$ denotes the positive eigenspan of $A_{b_t}$.  The
following lemma shows
that such a  path can always be found, and that the resulting spectral flow is
independent of the choice of the Riemannian metric. What makes the
proof of Lemma
\ref{apslem} tricky is that we do not assume the kernels of the family
$A_{b_t}$   have constant dimension.

\begin{lemma}\label{apslem} Suppose that $B_0$ and $B_1$ are two flat $U(n)$
connections on $X$ whose holonomies $\al_0,\al_1$ lie in the same path
component of
$\chi(\pi_1X, U(n))$.

Then, perhaps after gauge transforming $B_1$,  there is a continuous piecewise
smooth  path
$B_t$ of flat
$U(n)$ connections joining them
and a  continuous piecewise smooth path   $P(t)=F^+_{b_t}(t)\oplus W(t)$ of
self-adjoint APS boundary conditions for the corresponding odd signature
operators, with
$P(0)=F^+_{b_0}\oplus V_{X,\al_0}$ and $P(1)=F^+_{b_1}\oplus V_{X,\al_1}$.
In fact this path   can be taken to be piecewise real analytic.

\begin{sloppypar}
Moreover, the spectral flow of the path $D_{B_t,P(t)}$ of self-adjoint
operators, $\SF(D_{B_t,P(t)})\in \zz$, depends only on the path
$\al\co I\to \chi(\pi_1X,U(n))$ of holonomies of $B_t$ and the choice of
$W(t)$. In particular it is
  independent
of the choice of Riemannian metric  on $X$ (and $\del X$).
\end{sloppypar}
\end{lemma}

\begin{proof} We prove the last assertion first.  Note that
the kernel of $D_{B_t,P(t)}$   is isomorphic to
$$\text{image} \big(H^*(X,\del X; \cc^n_{\al_t})\to H^*(X; \cc^n_{\al_t})\big)
\oplus \big(W(t)\cap \ga(V_{X,\al_t})\big)$$
by \eqref{eq1.6} (as usual $\al_t$ denotes the holonomy representation of
$B_t$). This is not quite a homotopy invariant since it is possible that by
varying the metric one could change the intersection of $W(t)$ (which is metric
independent) with $\ga(V_{X,\al_t})$ (which can vary with the metric since
$\ga$ does). However, at the endpoints $W(i)=V_{X,\al_i}, i=0,1$, and since
$\ga(L)=L^\perp$ for any Lagrangian, $W(i)\cap \ga(V_{X,\al_i})=0$ for $i=0,1$.
Thus the dimension of the kernel of $D_{B_i,P(i)}$ is independent of the
choice of Riemannian metric on $X$ and $\del X$.

Varying the Riemannian metric on $X$ varies the path $D_{B_t,P(t)}$
continuously in the space of self-adjoint  operators (\cite[Sec.\ 3]{BooLesPhiUFS}).
Moreover, since the
dimensions of the kernels at the endpoints are homotopy invariants   it
follows that varying the metric does not change the spectral flow.

Suppose that $B_t'$ is another continuous path of flat connections with the
same holonomy as $B_t$.  Then one can find a continuous path of gauge
transformations $g_t$ so that $g_t\cdot B_t=B_t'$.  It is then straightforward
to use $g_t$ to define   unitary transformations  which conjugate the
operators $D_{B_t,P(t)}$ to $D_{B_t',P'(t)}$, and so the spectral flows of the
two paths coincide. We leave the details to the reader

We turn to the problem of constructing the paths $B_t$ and $W(t)$.
Using the main result of \cite{fine-kirk-klassen} one can
find a piecewise real analytic path $B_t$ of flat connections joining
$B_0$ to $B_1'$, where $B_1'$ is a flat connection gauge equivalent to $B_1$.
Indeed we can pick a piecewise analytic path $\al_t$ in $\chi(\pi_1X,U(n))$
from $\al_0$ to $\al_1$,  and the main theorem of \cite{fine-kirk-klassen}
shows that one can find a finite covering of the interval $[0,1]$ and analytic
paths of flat connections in each subinterval with the corresponding holonomy.

By relabeling assume that $B_1=B_1'$.  Then by subdividing further if
necessary, the interval
$[0,1]$ can broken down into subintervals
$[t_i,t_{i+1}]$ so that  on each such subinterval:
\begin{enumerate}
\item The path $B_t$ is real-analytic.
\item The kernel of the tangential operator, $\ker A_{b_t}$, has constant
dimension on the interior of the interval.
\end{enumerate}
The reason why the second condition can be met is that the subspaces
$$S_k:=\{ \al\in\chi(\pi_1X,U(n)) | \dim H^*(\del X,\cc^n_\al)\leq k\}\subset
\chi(\pi,U(n))$$
form a real-analytic subvariety  and so an analytic path intersects
it in a finite number of points.

For convenience, reparametrize the path so that on alternate intervals the path
is constant, i.e.\ $B_t$ is constant on $[0,t_1], [t_2,t_3], \cdots ,
[t_{2m},t_{2m+1}], \cdots , [t_{2\ell},1]$.

We now construct the path $W(t)$.

  The $W(i)$ have already been chosen at the endpoints: we take
$W(0)=V_{X,\al_0}\subset
  H^*(\Si;\cc^n_{\al_0})$ and $
W(1)=V_{X,\al_1}\subset  H^*(\Si;\cc^n_{\al_1}).$

Next, on each interval $[t_{2m-1},t_{2m}]$,  by the Kato selection lemma
\cite{kato} we can find analytically varying eigenvectors
$\psi_j(t)$ with analytically varying eigenvalues $\mu_j(t)$ for $j\in
\{\pm1,\pm2,\cdots \}$.

The set $\{\psi_j\}$ can be partitioned into the finite subset
$$K:= \{\psi_j\ | \ \mu_j(t)=0 \text{ for all } t\in [t_{2m-1},t_{2m}]\}$$
and its complement. Moreover, by relabeling we may assume that
$$K=\{ \psi_j\ | \ j=\pm1, \pm 2, \cdots, \pm \ell\}.$$
We   further   assume, by changing bases, that
  $\ga(\psi_j(t))=\psi_{-j}(t)$ for $j=1,\cdots ,\ell$ and $\psi_j(t)\perp
\psi_k(t)$ for
$ j,k\in \{\pm 1,\cdots \pm\ell\}$ for all $t\in [t_{2m-1},t_{2m}]$.  This is
possible because $$S(t):=\text{span}\{ \psi_j(t)\ |\ j=\pm1,\cdots\pm\ell\}$$
is preserved by $\ga$ for each $t\in [t_{2m-1},t_{2m}]$ and hence is an
analytically varying family of Hermitian symplectic  spaces   (note
that
$S(t)=\ker A_{b_t}$ for $t$ in the interior of $[t_{2m-1},t_{2m}]$).
The spaces $S(t)$ contain an analytically varing family of Lagrangian
subspaces
$$L(t):=\text{span}\{ \psi_j(t)\  | \ j=1,\cdots, \ell\}\subset S(t).$$
Then define $W(t)$ on $[t_{2m-1},t_{2m}]$ as follows.
\begin{enumerate}
\item For $t_{2m-1}<t<t_{2m}$ take
   $$W(t)=L(t).$$
\item  For $t=t_{2m-1}$,
\[\begin{split}
&W(t_{2m-1})= L(t_{2m-1}) \oplus\\ 
&\oplus\text{span}\big\{ \psi_j(t_{2m-1})\,\big|\,
\mu_j(t_{2m-1})=0 \text{ and } \mu_j(t)> 0 \text{ for }  t_{2m-1}< t<t_{2m}\big\}.
  \end{split}
\]
Briefly, $W(t_{2m-1})$ is the  span of $\lim_{t\to t_{2m-1}^+}L(t)$ and those
zero eigenvectors that deform to positive  eigenvectors for $t>t_{2m-1}$.
\item For $t=t_{2m}$,
\[\begin{split}
&W(t_{2m})= L(t_{2m}) \oplus \\
&\oplus\text{span}\big\{ \psi_j(t_{2m})\,\big|\,
\mu_j(t_{2m})=0
\text{ and } \mu_j(t)> 0 \text{ for }  t_{2m-1}< t<t_{2m}\big\}.
  \end{split}
\]
Thus $W(t_{2m})$ is the  span of $\lim_{t\to t_{2m-1}^-}L(t)$ and those zero
eigenvectors that deform to positive  eigenvectors for $t<t_{2m}$.
\end{enumerate}

By construction, $P(t)=F_{b_t}^+(t)\oplus W(t)$ is smooth (even analytic) on
the  interval $[t_{2m-1},t_{2m}]$. That the $W(t)$ are Lagrangian is immediate
except possibly at the endpoints. But at the endpoint  $t_{2m-1}$ the
decomposition of
$\ker A_{b_{t_{2m-1}}}$ into the nullvectors that ``stay null'' and those that
deform into non-zero eigenvectors is a decomposition as a  symplectic direct
sum, and in the summand corresponding to the eigenvectors that deform into
non-zero eigenvectors the subspace of those that deform into positive
eigenvectors is Lagrangian. Thus $W(t_{2m-1})$ is a direct sum of Lagrangian
subspaces, and hence is Lagrangian. A similar argument applies to
$W(t_{2m})$.

It remains to define $W(t)$ on the intervals $[t_{2m}, t_{2m+1}]$. On these
intervals the connection $B_t$ is constant. Hence the symplectic space $\ker
A_{b_t}$ is constant also.  The Lagrangians $W(t_{2m})$ and $W(t_{2m+1})$ have
already been defined, so just pick some smooth path $W(t)$
interpolating between
these two.  With care this path can be chosen to be analytic.  Clearly the
path $P(t)=F_{b_t}^+(t)\oplus W(t)$ is smooth on this interval.
\end{proof}

The following theorem  says that  the spectral flow, which by Lemma
\ref{apslem} is a diffeomorphism invariant,   is in fact an invariant of
homotopy equivalences  which restrict to diffeomorphisms on the boundary.

\begin{thm}\label{sfthm} Suppose that $F\co X'\to X$ is a homotopy equivalence
which restricts to a diffeomorphism $f=F|_{\del X'}\co\del X'\to \del X$.  Assume
that $B_t$ is a continuous, piecewise smooth path of flat $U(n)$  connections
on $E\to X$. Use $F$ to pull back the path $B_t$ to a path of flat  connections
$B'_t$ on $X'$ and to identify $\del X$ with $\del X'$, and choose a path
$P(t)$ of APS boundary conditions as in Lemma \ref{apslem}.

Then
$$\SF(X,D_{B_t,P(t)})_{t\in [0,1]}=  \SF(X',D_{B'_t,P(t)})_{t\in [0,1]}.$$
\end{thm}
\begin{proof}
Since spectral flow is additive with respect to composition of paths, we may
assume that the path   $D_{B_t,P(t)}$ is a smooth path of self-adjoint
operators. Denote this path  by $D_t$;   thus
$\eta(D_t)=\eta( D_{B_t, W(t)},X)$.

Theorem \ref{thm1.1}, part 1 and the definitions imply that
\begin{equation}\label{eq7.6}
\begin{split}
\eta(D_s)- \eta(D_0) 
&= \eta( D_{B_s, V_{X,\al_s}}, X)+m(\ga(V_{X,\al_s}),
W(s))_{(\al_s,g)}\\
&\quad -\eta( D_{B_0,V_{X,\al_0}},X)
- m(\ga(V_{X,\al_0}), W(0))_{(\al_0,g)} \\
   &=\rho(X, D_{B_s},g)+m(\ga(V_{X,\al_s}), W(s))_{(\al_s,g)}\\
&\quad -\rho(X,D_{B_0},g) - m(\ga(V_{X,\al_0}), W(0))_{(\al_0,g)},
\end{split}
\end{equation}
where
$\al_t$ denotes the holonomy of $B_t$ and $g$ is the metric on $\del X$.

The reduction of  the $\eta$--invariants $\eta(D_t)$
 modulo $\zz$ is smooth in $t$. Combining the
formula  (see e.g.~\cite[Lemma 3.4]{Kirk-Lesch})
\begin{equation}\label{eq7.10}
\eta (D_s)-\eta (D_0)= 2\SF(D_t)_{t\in [0,s]}
-(\dim\ker D_s-\dim\ker D_0)+
\int_0^s
\frac{d\eta(D_t)}{dt}  dt\end{equation}
   with \eqref{eq7.6} yields
\begin{equation}\label{eq7.7}
\begin{split}
\rho(X, &\al_s, g)- \rho(X,\al_0,g)\\
&\quad + m(\ga(V_{X,\al_s}), W(s))_{(\al_s,g)}-
m(\ga(V_{X,\al_0}), W(0))_{(\al_0,g)} \\
& =2\SF(D_t)_{t\in[0,s]}-(\dim\ker D_s-\dim\ker
D_0) +\int_0^s\frac{d\eta(D_t)}{dt} dt.
\end{split}
\end{equation}
Similarly
\begin{equation}\label{eq7.8}
\begin{split}
\rho(X', &\al_s, g)- \rho(X',\al_0,g)\\
&\quad + m(\ga(V_{X',\al_s}), W(s))_{(\al_s,g)}-
m(\ga(V_{X',\al_0}), W(0))_{(\al_0,g)}\\
& =2\SF(D'_t)_{t\in[0,s]} -(\dim\ker D'_s-\dim\ker
D'_0)+\int_0^s
\frac{d\eta(D_t')}{dt} dt.
\end{split}
\end{equation}
where $D'_t$ denotes the odd signature operator on $X'$ coupled to the path
$B'_t$ with boundary conditions given by the projection to $P(t)$.

    Taking the  difference of \eqref{eq7.7} and
\eqref{eq7.8} and using  Theorem
\ref{FLWKL}, the fact that
$V_{X,\al_s}=V_{X',\al_s}$, and \eqref{eq1.6} one concludes
$$ 2\SF(D_t)_{t\in[0,s]} -
2\SF(D_t')_{t\in[0,s]}=\int_0^s \frac{d\eta(D_t')}{dt} dt -
\int_0^s \frac{d\eta(D_t)}{dt} dt.$$
The left side is an integer-valued function of $s$. The right side is a smooth
real-valued function of $s$ which vanishes at $s=0$. Thus both sides vanish for
all $s$ and so
$$ \SF(D_t)_{t\in[0,s]}=
   \SF(D_t')_{t\in[0,s]}$$
as  desired.
\end{proof}

Theorems \ref{FLWKL} and \ref{sfthm} do not hold without some assumption about
the restriction of the  homotopy equivalence to  the
boundary. Here is an an example. Consider the complements
$X=S^3-$nbd$(K)$ of the Square knot $K$  and $X'=S^3-$nbd$(K')$ of the Granny
knot $K'$. (The Square knot is a connected sum of a right-handed Trefoil knot
and a left-handed Trefoil knot. The Granny knot is the connected
sum of two right-handed Trefoil knots.)  The
spaces $X$ and $X'$  have isomorphic fundamental groups and are
aspherical, and so they are homotopy equivalent (see e.g.~\cite{rolfsen} and \cite{hempel}). Each
has a torus boundary. But  there does not exist a homotopy equivalence
which restricts to a diffeomorphism (or even  a homotopy equivalence) on the
boundary.  This follows from Waldhausen's theorem \cite{waldhausen} (it also
follows from the following argument  and Theorem \ref{sfthm}).

Since $H_1(X;\zz)=\zz=H_1(X';\zz)$, the $U(1)$ character variety
of $\pi_1X=\pi_1X'$ is a circle, parameterized by the image
$z=e^{ix}$ of the generator of  the first homology.  Fix a generator
$\mu\in H_1(X;\zz)$. Let $F\co X\to X'$ be a
homotopy equivalence.
   Let
$\al_z\co\pi_1(X)\to U(1)$ be the representation which takes $\mu$
to $z=e^{ix}$. Let $B_z$ be a path of flat connections on $X$ with holonomy
$\al_z$ and restriction $b_z$ to the boundary.  Let $B_z'=F^*(B_z)$.

\begin{thm}\label{squaregranny} The 3-manifolds $X$ and $X'$ are homotopy
equivalent, have diffeomorphic boundaries, but for no choice of metrics
$g$,
$g'$ on $\del X$, $\del X'$  is
$\rho(X,\al ,g)$ equal to $\rho(X',\al ,g')$ for all $\al \in
\chi(\pi_1X,U(1))$.
\end{thm}
\begin{proof}
Suppose that $g$ is a metric on $\del X$ and $g'$ a
metric on $\del
X'$.

   Lemma \ref{lem3.1} shows that for $z\ne 1$,
$H^*(\del X;\cc_{\al_z})=\ker A_{b_z}=0$.
Thus the boundary conditions given by the projection to the positive
eigenspan $F^+_{b_z,g}$ of $A_{b_z}$ are self-adjoint (for $z\ne 1$) and vary
smoothly in $z$.  Similarly for
$X'$.

The kernel of the operator $D_{B_z}$ with these boundary
conditions is zero except for those $z$ which are roots of the
Alexander polynomial of the Granny knot $K$. Moreover, the
spectral flow as $z$ moves through the root is given by the change
in the Levine-Tristram signature (see e.g.~\cite{chuck}) of $K$.
These facts are proven e.g.\ in \cite{KKR}.

   Since the  Square and Granny  knots have
isomorphic groups, their Alexander polynomials are the same, namely
$(z^2-z+1)^2$. Thus eigenvalues of $D_{B_z}$ and $D_{B'_z}$
cross zero    for the same values of $z$, namely $e^{2\pi i/6}$ and $e^{5\pi
i/6}$.

  The spectral flow  through these values of $z$ is
different for $X$ and $X'$. Indeed,   the Square knot is slice and
hence has vanishing Levine-Tristram signatures; thus
$\SF(D_{B_z})$ is zero.  But the Levine-Tristram signatures for
the Granny knot are non-trivial (they detect the non-sliceness of
the Granny knot). Thus if $J$ is a small interval in $U(1)$
containing a root of the Alexander polynomial (to be explicit,  we
can take $J=\exp(2\pi
i[\tfrac{1}{6}-\epsilon,\tfrac{1}{6}+\epsilon])$)
$$\SF(X, D_{B_z,  F^+_{b_z}})_{z\in J}=0 \text{ and } \SF(X, D_{B_z ,
F^+_{b'_z}})_{z\in J}= 2.$$
   The reduction of $\rho(X,\al_z,g)$ to $\rr/\zz$ is
continuous in $z$. Since the integer jumps of $\rho(X,\al_z,g)$ as $z$ varies
are given by the spectral flow (see
\eqref{eq7.10}), it follows that for some $z$ near  $e^{2\pi i/6}$,
 \[\rho(X,\al_z,g)\ne\rho(X',\al_z,g').\tag*{\qed}\]
\exendproof

\section{Determinant bundles and variation of the\newline
   $\rho$--invariant mod $\zz$ on manifolds with boundary}
\label{det bundles}

As before let $X$ be a compact odd-dimensional manifold with boundary $\Si$.
We are given a bundle $E\to X$ and a flat connection $B$ on $E$ in cylindrical
   form $b + \tfrac{\del}{\del u}$  near the boundary, and a Riemannian
   metric $\tilde{g}$ on $X$  in cylindrical form $g + du^2$ near the  boundary,
   for some Riemannian metric $g$ on $\Si$.

Recall that the {\em determinant line} of the operator $A_b$ is the
   (complex) vector
space
$$\det(\ker A_b)= (\det\ker A_b^+)^{-1}\otimes \det\ker A_b^-,$$
where $\ker A_b^{\pm}:=\ker A_b\cap\ker(\gamma\mp i)$, i.e.\
$\ker A_b$ is considered as a $\zz_2$--graded vector space with
grading operator $-i\gamma$. For details about graded determinant
lines we refer to \cite[Sec.\ II]{Dai-Freed}.

Given a Lagrangian $W\subset \ker A_b$, it can be written as the graph of
the unitary isomorphism $\phi(W)\co\ker A_b^+\to\ker A_b^-$. Thus
$\det(\phi(W))$ is naturally an element of  $\det(\ker A_b)$. Theorem
\ref{thm1.1}
shows that
\begin{equation}\label{G2.96}
e^{2\pi i \etab(D_{B,W},X)}\det(\phi(W))^{-1}=
e^{2\pi i \etab(D_{B,V_{X,\al}},X)}\det(\phi(V_{X,\al}))^{-1}.
\end{equation}
Hence the expression
$s(X,\alpha,g):=e^{2\pi i \etab(D_{B,W},X)}\det(\phi(W))^{-1}$ is independent
of $W$, and so well--defines an element  in $(\det \ker A_b)^{-1}$.
As a consequence,
we obtain a well--defined element
\begin{equation}\label{G2.97}\begin{split}
&e^{2\pi i \etab(D_{B,W_\al},X)}
e^{-2\pi i\etab(D_{\Theta,W_\tau},X)}
\det(\phi(W_\al))^{-1}\otimes \det(\phi(W_\tau))\\
&\in (\det \ker A_b)^{-1}\otimes \det \ker A_\theta,
			     \end{split}
\end{equation}
independent of the choice of $W_\al$ and $W_\tau$.
Here $\Theta$ denotes the trivial connection on $\cc^n\times X\to X$ and
$\theta$ its restriction to $\Si$. By slight abuse of notation we will denote
the element
$s(X,\alpha,g)s(X,\tau,g)^{-1}$ of
   $ (\det \ker A_b)^{-1}\otimes \det \ker A_\theta$ given in  \eqref{G2.97}
by $e^{\pi i \rho(X,\al,g)}$.

For the spin Dirac operator the fact that the expression
on the left hand side of \eqref{G2.96} gives a well--defined element
of the inverse determinant line was observed first by X. Dai and
D. Freed \cite[Sec.\ I]{Dai-Freed}. As pointed out in \cite{Dai-Freed}
the result easily transfers to general Dirac operators. Due to
the exponentiation, \eqref{G2.96} does not need the full strength
of Theorem \ref{thm1.1}. As in \cite{Dai-Freed} it can
be derived already from \cite{LesWoj:IGA} where the dependence
of the $\operatorname{mod}\, \zz$ reduced $\eta$--invariant on the
boundary condition
is investigated.

We now follow Dai and Freed in \cite{Dai-Freed} and generalize
\eqref{G2.96} and \eqref{G2.97} to the parameterized context.
As a general reference for the analysis of elliptic families
we refer to the book \cite{BGV}.

As parameter space we take
$$\mathscr{P}_\Si= \mathscr{F}_\Si\times\mathscr{M}_\Si, $$
where
$\mathscr{M}_\Si$ denotes the space of Riemannian metrics on $\Si$ and
$\mathscr{F}_\Si$ denotes the space of flat connections on the bundle
$E|_\Si\to \Si$.

Given a manifold $X$ with $\del X=\Si$, we will also use the parameter space
$$\mathscr{P}_X=\mathscr{F}_X\times\mathscr{M}_\Si,$$
where $\mathscr{F}_X$ denotes the space of flat connections on $E\to X$
in cylindrical form $B=b +\tfrac{\del}{\del x}$ on a collar.
Notice that we
take the space of Riemannian metrics on $\Si$, not $X$, to define
$\mathscr{P}_X$.

Also note that in contrast to \cite{Dai-Freed} the
parameter spaces $\mathscr{P}_\Si$ and $ \mathscr{P}_X$ are
infinite--dimensional. For the discussion of the determinant line
bundle, however, this does not cause any additional difficulties.
The reader who prefers not to worry about the manifold structure
of $\mathscr{P}_X, \mathscr{P}_\Sigma$   may think of having chosen a finite--dimensional
submanifold of $\mathscr{P}_X, \mathscr{P}_\Sigma$.

To be more precise we consider the trivial fibration
\begin{equation}
\pi\co X\times \mathscr{P}_X\to \mathscr{P}_X.
\end{equation}
Give the fiber $X_{p}$ over $p=(B,g)$ a Riemannian metric so that in a
fixed collar of $\partial X$ the metric takes the form
  $g +du^2$. This can be done in a smooth way over $\mathscr{P}_X$; for
example, fix a metric on the interior of $X$ and use a cutoff
function in a slightly larger collar to interpolate between this fixed
metric and $g+ du^2$. Call the resulting extended metric $\tilde{g}$.
The fiber of the relative tangent
bundle $T(X\times \mathscr{P}_X/\mathscr{P}_X):=\ker T\pi\co
T(X\times \mathscr{P}_X)\to T\mathscr{P}_X$ over the point $(x,(B,g))
\in X\times\mathscr{P}_X$ is given by $T_{x}X\oplus 0$ and the
metric on $T_{(x,(B,g))}(X\times \mathscr{P}_X/\mathscr{P}_X)$
is $\tilde g_x$.

Since $\pi$ is a product we have a natural horizontal structure which
is given by the kernel of the tangent map of the projection
$X\times\mathscr{P}_X\to X$ onto the first factor.

Summing up we have a Riemannian structure on the fibration $\pi$
in the sense of \cite[p.\ 5159]{Dai-Freed} resp.
\cite[Sec.\ 10.1]{BGV}.

\begin{sloppypar}
Given $p=(b,g)\in \mathscr{P}_\Si$ we modify the previous
   notation slightly and denote by
   $A_p$ the odd signature operator coupled to $b$ in order to emphasize its
   dependence on the Riemannian metric $g$. Then $A_p$ acts on the bundle
$(\Lambda\mathscr{T})_p:=\Lambda T^*(\Sigma;E):=\oplus_{q\ge 0}\Lambda^qT^*(\Sigma;E)$,
where the notation $(\Lambda\mathscr{T})_p$  also emphasizes the dependence on
$p$ through the metric. In fact $(\Lambda\mathscr{T})_p$ is the restriction of
the bundle
$\Lambda\mathscr{T}:=\Lambda T^*(\Sigma\times \mathscr{P}_\Sigma/\mathscr{P}_\Sigma,E):=
\oplus_{q\ge 0}\Lambda^qT^*(\Sigma\times \mathscr{P}_\Sigma/\mathscr{P}_\Sigma,E)$
(the exterior power bundle of the relative cotangent bundle)
to the  fiber  $\Sigma\times \{p\}$. Thus, the
 $(A_p)_{p\in\mathscr{P}_\Sigma}$ form a smooth family of Dirac type operators
in the sense of \cite[Sec.\ 9.2]{BGV}.
\end{sloppypar}

We   then form the associated determinant bundle $\det \ker A\to
\mathscr{P}_\Si$
   whose fiber
over $p\in \mathscr{P}_\Si$ is
$\det \ker A_p=(\det \ker A_p^+)^{-1}\otimes \det\ker A_p^-$.
   One can find a
   careful construction of  $\det \ker A$ in many articles, e.g.\
   \cite{BGV, Bismut-Freed}, as well as the construction of the
Quillen metric and connection on $\det \ker A$.  We outline briefly
   the reason why the fibers $\det \ker A_p$ glue together to form a
smooth bundle
for the benefit of the reader. The space $\mathscr{P}_\Si$ is
   covered by open sets $U_\mu$,
   $\mu>0$,\label{twentyfive}
consisting of those $p$ so that
   $\mu\not\in \Spec (A_p)$.  If $F_p^\pm(\mu)$ denotes the span of those
   eigenvectors of $A_p$ whose eigenvalues lie in $(-\mu,\mu)\setminus\{0\}$,
then the vector
spaces $H(\mu)_p$ defined by
$$H(\mu)_p:=F_p^+(\mu)\oplus \ker A_p\oplus F_p^-(\mu)$$
   form  a smooth, finite dimensional
   vector bundle over $U_\mu$, whose fibers are invariant under $\gamma_p$,
   since $\gamma_p A_p=-A_p\gamma_p$. Thus $H(\mu)_p$
is a Hermitian symplectic space and  so has a decomposition
\begin{equation}
H(\mu)_p =\bigl(H^+(\mu)_p\oplus\ker A_p^+\bigr)\oplus\bigl(\ker
A_p^-\oplus H^-(\mu)_p\bigr)
\label{G8.4}
\end{equation}
into the $\pm i$ eigenspaces of $\gamma_p$ acting on $H(\mu)_p$.
The decomposition \eqref{G8.4} yields
\begin{equation}
   \det H(\mu)_p=\det(\ker A_p)\otimes \det(H^+(\mu)_p)^{-1}\otimes
    \det(H^-(\mu)_p).
\end{equation}
The crucial observation now is that $\det(H^+(\mu)_p)^{-1}\otimes
    \det(H^-(\mu)_p)$ is \emph{canonically} trivial since $\det(A_p^+\co H^+(\mu)_p
\to H^-(\mu)_p)$ is a canonical nonzero section of
$\det(H^+(\mu)_p)^{-1}\otimes\det(H^-(\mu)_p)$. Consequently
$\det H(\mu)_p$ is canonically isomorphic to $\det \ker A_p$. Identifying
$\det \ker A|U_\mu$ with $\det H(\mu)|U_\mu$ shows that $\det \ker A$ indeed
is a smooth line bundle over $\mathscr{P}_\Sigma$.

For future reference we denote the canonical (bundle) isomorphism
\begin{equation}
       \varphi(\mu)\co\det(\ker A)|U_\mu\longrightarrow \det H(\mu)|U_\mu,\quad
            \xi_p\mapsto \xi_p\otimes \det(A_p^+|H^+(\mu)).
\label{G8.6}
\end{equation}
From \eqref{G2.96} one now   infers (cf.\ \cite{Dai-Freed}):

\begin{prop}\label{DFprop} Given $p=(b_p,g_p)\in\mathscr{P}_\Si$, let
$A_{\theta,p}$ denote the odd signature operator on the trivial bundle
$\cc^n\times \Si\to
\Si$ with respect to the Riemannian metric $g_p$ on $\Si$.

Then
the
vector spaces
$\det(\ker A_p)^{-1}\otimes \det(\ker A_{\theta,p})$   form a
smooth vector bundle
$$\det(\ker A)^{-1}\otimes \det(\ker A_\theta) \to \mathscr{P}_\Si.$$
   Moreover, if $\del X=\Si$ the $\rho$--invariant defines a smooth lift
\[\begin{diagram}\dgmag{900}
\node[2]{\det(\ker A)^{-1}\otimes \det(\ker A_\theta) }\arrow{s}\\
\node{\mathscr{P}_X}\arrow{e,t}{i^*}\arrow{ne,t}{e^{\pi i \rho}}
\node{\mathscr{P}_\Si}
\end{diagram}\]
where $i^*$ denotes the restriction.
\end{prop}
\begin{proof} From the previous discussion and in view of \eqref{G2.96}
the statement of the proposition is clear except the fact
that $e^{\pi i\rho}$ defines a \emph{smooth} lift into
$\det(\ker A)^{-1}\otimes\det(\ker A_\theta)$.
Although this fact was proved in
\cite{Dai-Freed} it also follows from our results. The key is that
\eqref{G2.96} can be generalized in such a way that one obtains
smooth sections over $U_\mu$. For a Lagrangian $W\subset H(\mu)_p$ denote
by $\eta(D_{B,W},X)$ the $\eta$--invariant of $D_B$ with respect to
the boundary condition given by the orthogonal projection onto
$(F_p^+\ominus F^+_p(\mu))\oplus W$. Denote by $\phi(\mu)_p(W)$ the analogue
of $\phi(W)$ for the space $H(\mu)_p$. Then $\det(\phi(\mu)_p(W))$ is a
canonical element of $\det H(\mu)_p$.

Given two such subspaces $W_1, W_2$ then by \cite[Thm.\ 4.2]{Kirk-Lesch}
we have the following generalization of \eqref{G2.96}:
  \begin{equation}\label{G2.96a}
e^{2\pi i \etab(D_{B,W_1},X)}\det(\phi(\mu)_p(W_1))^{-1}=
e^{2\pi i \etab(D_{B,W_2},X)}\det(\phi(\mu)_p(W_2))^{-1}.
\end{equation}
This shows that the expression $e^{2\pi i
\etab(D_{B,W_1},X)}\det(\phi(\mu)_p(W_1))^{-1}$ is independent of $W_1$ and
choosing a smooth family of Lagrangians in $(H(\mu)_p)_{p\in U_\mu}$ over
$U_\mu$ gives a smooth section of $\det(\ker A)$.\end{proof}

The (tensor product of two) determinant
bundle(s)  
\[\det(\ker A)^{-1} \otimes \det(\ker A_\theta)\to \mathscr{P}_\Si\]
admits the Quillen metric \cite{Quillen} and its natural compatible connection
$\nabla^Q$ \cite{Bismut-Freed}. The main result of \cite{Dai-Freed} can be used to
compute  $\nabla^Q(e^{i\pi \rho})$.

\begin{sloppypar}
Let us first briefly recall the main facts about metrics and connections
on $\det(\ker A)$ (resp. $\det(\ker A)^{-1}\otimes \det(\ker A_\theta))$.
We use the
  notation  from page \pageref{twentyfive}. Since the fiber of the
relative tangent bundle $T(\Sigma\times\mathscr{P}_\Sigma/\mathscr{P}_\Sigma,E)$
over $(x,p)\in \Si\times \mathscr{P}_\Sigma$ is $T_x\Si\oplus 0$, the relative
tangent bundle is naturally a Riemannian vector bundle. Consequently,
$\Lambda\cT$ inherits a natural metric from the relative tangent bundle.
Furthermore, by \cite[Prop.\ 10.2]{BGV} $T(\Sigma\times\mathscr{P}_\Sigma/\mathscr{P}_\Sigma,E)$
has a natural connection which is induced solely by the Riemannian structure of
the fibration $\pi$. This connection induces a
connection on the bundle $\Lambda\cT$. Finally, we note that for each $p$
the metric on $\Lambda\cT$ induces an $L^2$--structure on sections
of $(\Lambda\cT)_p$. In the terminology of family index theory
the space $C^\infty(\Sigma\times\{p\},\Lambda\cT_p)$ of sections
of $(\Lambda\cT)_p$ is viewed as the fiber over $p$ of an infinite--dimensional
Hermitian bundle $\pi_*\Lambda\cT$ whose sections are
$C^\infty(\Si\times \mathscr{P}_\Sigma,\Lambda\cT)$.
By \cite[Prop.\ 9.13]{BGV} there is a natural connection  $\nabla^{\pi_*\Lambda\cT}$
on $\pi_*\Lambda\cT$ which is compatible with the inner product. For details
we refer to \cite[Chap.\ 9]{BGV}.
\end{sloppypar}

The bundle $H(\mu)$ is a finite--dimensional sub--bundle of
$(C^\infty(\Sigma\times\{p\},\Lambda\cT_p))_{p\in U_\mu}$ and $\nabla^{\pi_*\Lambda\cT}$
projects to a connection $\nabla^\mu$ on $H(\mu)$. Furthermore, $H(\mu)$ inherits
a metric $h$ from the $L^2$--structure on
$(C^\infty(\Sigma\times\{p\},\Lambda\cT_p))_{p\in U_\mu}$.

However, the bundle map
$\varphi(\mu)$
\eqref{G8.6} is not an isometry with respect to the metric
$h$ and thus $h$ does in general not descend to a smooth metric on $\det\ker
A|U_\mu$. The norm of $\varphi(\mu)$ (with respect to the natural metrics on $\det\ker A|U_\mu$ and
on
$\det H(\mu)$) is given by
\[\begin{split}
\|\varphi(\mu)_p\|&=\bigl(\det A_p^+A^-_p|H^+(\mu)_p\bigr)^{1/2}\\
    &=\prod_{\lambda\in\Spec(A_p), 0<\lambda<\mu} \lambda.
   \end{split}
\]
Therefore, the \emph{Quillen metric} defined by
\[
     h^Q(\xi):=h(\xi) \prod_{\lambda\in\Spec A_p}\lambda:=
     h(\xi) \det\nolimits_\zeta(A_p^+A_p^-|\ker A_p^\perp)
\]
for $\xi\in\det\ker A_p$ is a smooth metric on $\det\ker A$. Here
$\det_\zeta(A_p^+A_p^-|\ker A_p^\perp)$ denotes the $\zeta$--regularized
determinant of the operator $A_p^+A_p^-|\ker A_p^\perp$.

The natural connection $\nabla^{\det,\mu}$ on $\det\ker A|U_\mu$ induced by $\nabla^\mu$
is in general not compatible with the Quillen metric.
  However,
there is a connection, $\nabla^Q$, on $\det\ker A$ which is compatible
with the Quillen metric. Formally, one has over $U_\mu$
\begin{align}  \nabla^Q&=\nabla^{\det,\mu}+\tr\bigl((A^+)^{-1}\nabla^\mu A^+\bigr)
\label{G8.7}\\
     &=:\nabla^{\det,\mu}+\beta_+^\mu.\nonumber
\end{align}
The right hand side of \eqref{G8.7} has to be suitably regularized;
for details we refer to \cite[Sec.\ 9.7]{BGV}.

As explained on page \pageref{twentyfive} the fibration $X\times \mathscr{P}_X$
has naturally the structure of a Riemannian fibration in the sense of
\cite[p.\ 5159]{Dai-Freed} resp. \cite[Sec.\ 10.1]{BGV}. Thus by
\cite[Prop.\ 10.2]{BGV} the relative tangent bundle
$T(X\times \mathscr{P}_X/\mathscr{P}_X)$ has a natural connection,
$\nabla^{T(X\times \mathscr{P}_X/\mathscr{P}_X)}$.
Let $R^X\in \Omega^2_{X\times \mathscr{P}_X}(\End{T(X\times\mathscr{P}_X/\mathscr{P}_X}))$ be the
curvature of this connection.

Next, let $\cE\to X\times \cP_X$ be the pullback of the bundle $E\to
X$ via the projection $X\times \cP_X\to X$. Let
   $\nabla^{\bf B}$  denote the connection on $\cE$ whose restriction to each
fiber $X_{B,g}$ is $B$. This can be constructed  by choosing an arbitrary
connection on $\cE$ and then adjusting it by the appropriate $1$-form.
The curvature of $\nabla^{\bf B}$,
   $F^{\bf B}\in \Omega^2_{X\times\cP_X}(\End \cE)$
restricts to zero in each fiber $X_{B,g}$ since $B$ is flat.

Similarly we construct the trivial connection $\nabla^{\bf \Theta}$ by
replacing $\cE$ by the trivial bundle and $B$ by the trivial
connection in the above formula.  Its curvature $F^{\bf \Theta}$
is zero and so the Chern character $\chern(F^{\bf \Theta})=n$ (as a form).

Then with these preparations,
the Dai-Freed theorem implies
the following.

\begin{prop}\label{prop6.1}
\begin{equation}
\nabla^Q(e^{i\pi \rho})= 2^{(\dim X-1)/2}\big( \int_X L(R^X)
(\chern(F^{\bf B})-
\chern(F^{\bf \Theta}))\big)_{[1]} \cdot e^{i\pi \rho},
\label{G8.10}
\end{equation}
where $\int_X$ denotes integration over the fibers and the subscript
``$[1]$'' means the 1-form part of an inhomogeneous differential form.
Furthermore,
$$L(R^X)=\det\nolimits^{1/2}\bigl(\frac{R^X/2}{\tanh(R^X/2)}\bigr)$$
denotes the Hirzebruch $L$--form associated to $R^X$.
\end{prop}
\begin{proof} In \cite{Dai-Freed} the theorem was stated for a
smooth family of spin Dirac operators. However, as they pointed out,
their result remains true for twisted Dirac operators. Since the
integrand in the right hand side of \eqref{G8.10}
is local and since every manifold is locally spin let us
assume for the moment that $X$ is spin. The complex Clifford algebra
$\cc l_{2k+1}$ (remember $\dim X=2k+1$) has two inequivalent irreducible
respresentations $\Delta^{\pm}$ and hence
\begin{equation}
\cc l_{2k+1}\cong \End(\Delta^+)\oplus \End(\Delta^-).\label{G8.9}
\end{equation}
Denote by $S_\cc^\pm\to X \times \mathscr{P}_X$ the spinor bundles
corresponding to
$\Delta^{\pm}$ associated to the relative tangent bundle $T(X\times
\mathscr{P}_X/\mathscr{P}_X)$. These inherit natural connections,
$\nabla^\pm$, from the connection
$\nabla^{T(X\times \mathscr{P}_X/\mathscr{P}_X)}$
on the relative tangent bundle
$T(X\times \mathscr{P}_X/\mathscr{P}_X)$.

From
the decomposition \eqref{G8.9} one easily infers (cf.\ also \cite[Sec.\ 4.1]{BGV}
for the even dimensional case) that the odd signature operator $D_\Theta$
is the spin Dirac operator coupled to the twisting bundle
$(S^+_\cc\otimes \cc^n,
\nabla^+\otimes\operatorname{id}+\operatorname{id}\otimes
\nabla^\Theta)$. Analogously,
$D_B$ is the spin Dirac operator coupled to the twisting bundle
$(S^+_\cc\otimes E,
\nabla^+\otimes\operatorname{id}+\operatorname{id}\otimes
\nabla^B)$.
Consequently, \cite[Theorem 1.9]{Dai-Freed} yields
\begin{equation}
\nabla^Q(e^{i\pi\rho})=\big( \int_X \hat
A(R^X)\wedge\chern(S_\cc^+,\nabla^+)\wedge
  (\chern(F^{\bf B})-
\chern(F^{\bf \Theta}))\big)_{[1]} \cdot e^{i\pi \rho}.
\end{equation}
It remains to identify the differential form $\hat
A(R^X)\wedge\chern(S_\cc^+,\nabla^+)$.
By the following Lemma we have
\begin{equation}
\hat A(R^X)\wedge\chern(S_\cc^+,\nabla^+)=2^k
L(R^X)=2^k\det\nolimits^{1/2}\bigl(\frac{R^X/2}{\tanh(R^X/2)}\bigr).
\label{G8.11}
\end{equation}
Note that although \eqref{G8.11} is an identity between differential forms,
it is in fact a statement about invariant polynomials on the special orthogonal
group and hence it follows indeed from the next lemma. \end{proof}

We point out that it is crucial in the following that $\hat
A(R^X)\wedge\chern(S_\cc^+,\nabla^+)=2^k L(R^X)$ is a Pontrjagin form and
hence is a sum of differential forms of degree divisible by four.

\begin{lemma} Let $M$ be a differentiable manifold and $E\to M$ a real oriented
vector bundle of rank $2k+1$ which carries a spin structure. Let $S_\cc^\pm(E)$
be the corresponding spinor bundles. Then
$\hat A(E)\wedge\chern(S_\cc^\pm(E))=2^k L(E).$
\end{lemma}
\begin{proof} A similar result for even rank bundles is well--known
(the lemma is probably  well--known also, however standard texts refer
to the even dimensional case only; see \cite[Sec.\ 4.1]{BGV}, \cite[Sec.\
3.3.5]{Gil:ITH}, \cite[Sec.\ III.11]{LawMic:SG}) and we will reduce the
lemma to the even dimensional case. By the splitting principle and since
the bundle is orientable we may assume that
\begin{equation}
   E\simeq \tilde E\oplus \rr_M,
\end{equation}
where $\tilde E$ is a real oriented bundle of rank $2k$ and $\rr_M$
denotes the trivial
$\rr$ bundle over $M$ (cf.\ \cite[Rem.\ III.11.3]{LawMic:SG}). Denote
by $S_\cc(\tilde E)$ the unique complex spinor bundle associated
to the spin structure on $\tilde E$. Then the representation theory
of the complex Clifford
algebras immediately implies that
\begin{equation}
   S_\cc^+(E)\simeq S_\cc^-(E)\simeq S_\cc(\tilde E)
\end{equation}
(isomorphisms as complex vector bundles). Hence we are reduced to the even
rank case and it follows (see the references above)
\[    \hat A(E)\wedge \chern(S_\cc^\pm(E))=\hat A(\tilde
E)\wedge\chern(S_\cc(\tilde E))=2^k L(\tilde E)=
    2^kL(E).\tag*{\qed}\]
\exendproof

Proposition \ref{prop6.1} implies the following.
\begin{thm}\label{thmcov} If one restricts to $SU(n)$ connections on
$E\to X$ or
if
$X$ is $(4\ell-1)$-dimensional  then
$$\nabla^Q(e^{\pi i \rho})=0.$$
\end{thm}

\begin{proof}
Recall that $X$ is a $(2k-1)$-dimensional manifold.
Decompose  the differential form $\chern(F^{\bf B})- \chern(F^{\bf
\Theta})$ into its
homogeneous
components:
$$\chern(F^{\bf B})- \chern(F^{\bf \Theta})=\chern_2(F^{\bf B}) +
\chern_4(F^{\bf B}) +
\chern_6(F^{\bf B})+\cdots$$
Similarly decompose
$$L(R^X)=L_0(R^X)+L_4(R^X) + L_8(R^X)+\cdots$$
Thus
\begin{equation}\label{eq6.65}\big(\int_X L(R^X)(\chern(F^{\bf B})-
\chern(F^{\bf
\Theta}))\big)_{[1]}=
\sum_{q\ge 1}\int_X L_{2k-2q}(R^X)\chern_{2q}(F^{\bf
B}). \end{equation}
Since the restriction of $F^{\bf B}$ to $X_p$ is flat,
$\chern_{2q}(F^{\bf B})=const(q)\text{Tr}((F^{\bf B})^q)$ has at
most
$q$ components in the ``$X$'' direction, i.e.\ writing $\chern_{2q}(F^{\bf
B})$ locally as a sum
   $$\chern_{2q}(F^{\bf B})_{x,p}= \sum_{i=0}^{2q}f^i_{x,p}
\alpha_i\wedge\beta_{2q-
i}$$
with $\alpha_i\in \Omega^i_X$ and $\beta_{2q-i}\in
\Omega^{2q-i}_{\cP_X}$, then $f^i=0$ for $i>q$.

This implies that the only possible non-zero summand in the right
side of \eqref{eq6.65} is the term with $q=1$, i.e.\
\begin{equation}\label{eq6.51}
\int_X L_{2k-2}(R^X)\wedge\chern_2(F^{\bf B}).
\end{equation}
But since $\chern_2(F^{\bf B})=c_1(F^{\bf B})$ and $L_{2k-2}(R^X)=0$
if $2k-2$ is
not divisible by
$4$, \eqref{eq6.51} vanishes if ${\bf B}$ is an $SU(n)$ connection or
if
$2k-2\ne 4\ell$. The result now follows from
Proposition \ref{prop6.1}.
\end{proof}

The following theorem exhibits a functoriality property of $\rho$ modulo $\zz$
for manifolds with boundary. It is closely related to   Theorem
\ref{FLWKL}, but the weaker hypothesis ($F$ need not be a homotopy equivalence)
gives a weaker conclusion: the $\rho$--invariants agree only modulo $\zz$.

\begin{thm}\label{cor16} Let $X$ and $X'$ be two odd dimensional manifolds and
suppose that $F\co X'\to X$ is a smooth map such that the restriction
$f=F|_{\del X'}\co\del X'\to \del X$ is a diffeomorphism. Let
$\al_0, \al_1\co\pi_1(X)\to SU(n)$ be two representations in the
same path component of $\chi(\pi_1(X),SU(n))$. Let $g_0$ and $g_1$
be two metrics  on $\del X$.
Then 
\[\begin{split}
    \rho(X,&\al_1,g_1)-\rho(X,\al_0,g_0)\\
  &\equiv 
\rho(X',F^*(\al_1),f^*(g_1))-\rho(X',F^*(\al_0),f^*(g_0))\pmod{\zz}.
  \end{split}
\]
In particular, if  $F\co X'\to X$ induces an isomorphism on
fundamental groups then  there is a factorization
\[\begin{diagram}\dgmag{430}
\node{\chi(\pi_1(X), SU(n))\times
\cM_\Sigma/\cD^0_\Si}\arrow[2]{e,t}{\rho(X)-\rho(X')}\arrow{se}\node[2]{\rr/\zz}
\\ \node[2]{\pi_0(\chi(\pi_1(X),
SU(n)))\times\cM_\Sigma/\cD^0_\Sigma}\arrow{ne,..}
\end{diagram}\]
and $\rho(X)-\rho(X')$ is zero on the path component of the
trivial representation. The result holds for $U(n)$ replacing
$SU(n)$ if $\dim X=4\ell-1$.
\end{thm}
\begin{proof} The map $F\co X'\to X$ induces a map $F^*$ on flat
connections $\mathscr{F}_X\to \mathscr{F}_{X'}$ by pulling back
connections. Using $f$ to pull back metrics on the boundary we see
that $F$ induces a map $F^*\co\cP_X\to \cP_{X'}$ so that
\[\begin{diagram}
\node{\cP_X}\arrow[2]{e,t}{F^*}\arrow{se,b}{i^*_X}
\node[2]{\cP_{X'}}\arrow{sw,b}{i_{X'}^*}\\ \node[2]{\cP_\Sigma}
\end{diagram}
\]
commutes.

Let $\al_t$ be a path of representations from $\al_0$ to $\al_1$.
Such a path can be chosen to be piecewise analytic and a
corresponding piecewise analytic path of flat connections $B_t$
with holonomy $\al_t$ can be found (\cite{fine-kirk-klassen}). By
adding the results in the end we may assume that the path $B_t$ is
analytic, and hence smooth. Let $g_t$ be a smooth path of metrics
on $\del X$ from $g_0$ to $g_1$. We identify $\del X$ and $\del
X'$ via $f$.

Choose a Lagrangian subspace $V_t$ in $\ker A_{(b_t,g_t)}$ for
each $t$ so that $V_0=V_{X,\al_0}$ and $V_1=V_{X,\al_1}$. (See Lemma
\ref{apslem}.) Similarly choose a Lagrangian subspace $W_t$ of $\ker
A_{(\theta,g_t)}$ with
$W_t=V_{X,\tau}$.

Let $c$ be the real number
$c=\rho(X,\al_0,g_0)-\rho(X',\al_0,g_0)$. Notice that by definition $c=0$ if
$\al_0$ is trivial. The smooth sections $$\phi_1\co t\mapsto\exp(\pi
i (\rho(X,\al_t, V_t,W_t),g_t))\det(\phi(V_t))\det(\phi(W_t))^{-1}$$
   and
   $$\phi_2\co t\mapsto
e^{c\pi i }\exp(\pi
i(\rho(X',\al_t,V_t),g_t))\det(\phi(V_t))\det(\phi(W_t))^{-1}$$ agree
at $t=0$ and by Theorem \ref{thmcov} satisfy
$\nabla^Q(\phi_1)=0=\nabla^Q(\phi_2)$ (since $e^{c\pi i }$ is
constant). In other words, $\phi_1$ and $\phi_2$ are two horizontal
lifts of the path $[0,1]\to \cP_\Sigma, t\mapsto (b_t,g)$. Since
they agree at $t=0$, they agree for all $t$. In particular, at
$t=1$ we conclude $$e^{\pi i \rho(X,\al_1,g_1)}=e^{\pi i
\rho(X',\al_1,g_1)}e^{c\pi i}.$$
This proves the first part of the theorem. The second part follows
from the first and the discussion following Theorem \ref{pii}.
\end{proof}
The second statement in Theorem \ref{cor16} should be compared to
\cite[Theorem 7.1]{Farber-Levine}.

We end this section with a discussion which shows that finding an
explicit dependence of the $\rho$--invariant on the metric on the
boundary is ultimately tied to the delicate construction of the
 connection $\nabla^Q$.

 Suppose that a representation
$\al\co\pi_1X\to U(n)$ is fixed and consider the function of metrics
on the boundary
\begin{equation}\label{roonm}\rho(X,\al)\co\cM_\Si\to \rr.\end{equation}
   This function is smooth, since $e^{i\pi\rho}$ is smooth and
since the dimension of the kernel of $D_{B,V_{X,\al}}$  is
independent of the metric by \eqref{eq1.6}.

\begin{prop}\label{matthiasremark} With the denotations of
\textnormal{\eqref{G8.7}} we have
\[\begin{split}
    d\rho(X,\alpha)&=-\frac{1}{\pi i} \frac{\det
\phi(V_{X,\alpha})}{\det\phi(V_{X,\tau})} \nabla^Q  \frac{\det
\phi(V_{X,\tau})}{\det\phi(V_{X,\alpha})}\\
     &=-\frac{1}{\pi i} \frac{\det
\phi(V_{X,\alpha})}{\det\phi(V_{X,\tau})} \nabla  \frac{\det
\phi(V_{X,\tau})}{\det\phi(V_{X,\alpha})}-\frac{1}{\pi i}\beta_+.
   \end{split}
\]
\end{prop}
\begin{proof} For the purpose of this proof the abusive notation $e^{i\pi\rho}$
introduced at the beginning of this section is too confusing. During this
proof we write $s(X,\alpha,g)s(X,\tau,g)^{-1}$
for the element of $(\det\ker A_b)^{-1}\otimes \det\ker A_\theta$
defined by \eqref{G2.97} and $e^{\pi i \rho}$ denotes the number
obtained by exponentiating
the $\rho$--invariant.

We cannot apply Theorem \ref{thmcov} directly. However,
since the parameter space is $\cM_\Si$ and $\alpha$ is fixed we have
$\chern(F^{\bf B})- \chern(F^{\bf \Theta})=0$ and hence by Proposition
\ref{prop6.1} and \eqref{eq6.65}
\begin{equation}
\nabla^Q(s(X,\alpha,g)s(X,\tau,g)^{-1})=0.
\label{G8.19}
\end{equation}
Over $\cM_\Si$ we have
\[  s(X,\alpha,g)s(X,\tau,g)^{-1}= e^{\pi i\rho(X,\alpha,g)}\det
\phi(V_{X,\alpha})^{-1}
\det\phi(V_{X,\tau}).\]
Consequently, \eqref{G8.19} implies
\begin{align*} d&\rho(X,\alpha,g)\\
&=-\frac{1}{\pi i}\det \phi(V_{X,\alpha})
\det\phi(V_{X,\tau})^{-1}\nabla^Q\bigl( \det \phi(V_{X,\alpha})^{-1}
\det\phi(V_{X,\tau})\bigr)\\
      &=-\frac{1}{\pi i}\det \phi(V_{X,\alpha})
\det\phi(V_{X,\tau})^{-1}\nabla\bigl( \det \phi(V_{X,\alpha})^{-1}
\det\phi(V_{X,\tau})\bigr)-\frac{1}{\pi i}\beta_+.\tag*{\qed}
   \end{align*}
\exendproof
Looking at the definition of $\beta_+$
we see that this result gives a link between the dependence of $\rho$ on the
metric and the variation of the (regularized) determinant of the
operator $A_b$.
This perhaps explains why we cannot expect the $\rho$--invariant to be
independent of the metric on the boundary.

\section{Topological  consequences}\label{Section 9}

It is known that the $\rho$--invariants
distinguish homotopy equivalent lens spaces \cite{wall-surgery}. By contrast, Neumann
showed in \cite{neumann} that $\rho$ is a homotopy invariant for
manifolds with free abelian fundamental
groups.

This leaves the problem of deciding exactly to what
extent the
$\rho$--invariant is a homotopy invariant open. One interesting  aspect of this
 problem is that it can be studied
one fundamental group at a time.

A conjecture of Weinberger
states (see
   \cite{weinberger}):

\newtheorem*{conja}{Conjecture A}

\begin{conja}[Weinberger]  If $M$ is a
closed
$(2k-1)$-manifold with torsion-free fundamental group then
$\rho(M,\al)$ depends only on the homotopy type of $M$.
\end{conja}

  Thus Neumann showed that Conjecture A holds for free abelian
groups. The Farber-Levine-Weinberger theorem solves the problem for
those groups  whose $U(n)$ character varieties are
connected, such as free groups. Wall's calculations for lens
spaces shows that the extension of the conjecture to all groups is
false: cyclic groups provide examples.

We make the following extension of the conjecture of Weinberger.
\newtheorem*{conjb}{Conjecture B}

\begin{conjb} Suppose that $F\co X\to X'$ is a
homotopy equivalence of manifolds with torsion free fundamental
groups  which
restricts to a diffeomorphism on the boundary. Endow
the boundaries with Riemannian metrics $g$, $g'$ so that the
restriction to the boundary is an isometry.

   Then for any  any
unitary representation $\al\co\pi_1X\to U(n)$
   $$\rho(X,\al,g)=\rho(X',\al,g').$$
\end{conjb}

This implies Weinberger's conjecture.
In this instance our Theorem \ref{FLWKL} implies Conjecture B for those
manifolds whose $U(n)$ character varieties are connected.

This reveals the following  strategy for attacking Conjecture A.

\begin{defn} We say a   homotopy equivalence $F\co M'\to M$ between closed
manifolds   {\em can be split along a separating hypersurface $\Si\subset M$} if, after a
homotopy of $F$,
\begin{enumerate}
\item $F$ is smooth and transverse to $\Si$ and the restriction
$F\co F^{-1}(\Si)\to \Si$   is a diffeomorphism, and
\item writing $M=X\cup_\Si Y$ and $M'=X'\cup_\Si Y'$,
$F$ restricts to homotopy equivalences $X'\to X$ and $Y'\to Y$.
\end{enumerate}

\end{defn}
The problem of determining when a homotopy equivalence can be
split along a hypersurface has been extensively studied; see 
e.g.~\cite[Chapter 12A]{wall-surgery} or \cite{cappell}.

We have the following result.

\begin{thm}\label{split} Let $F\co M\to M'$ be  a homotopy equivalence between
closed manifolds which can be split   along a hypersurface $\Si\subset
M$.  Write
$M=X\cup_\Si Y$.
Suppose that  $\al\co\pi_1M\to U(n)$ is a unitary representation  such
that the restriction  $\al|_X$ (resp.~$\al|_Y$) lies in the path component of
the trivial $U(n)$ representation of $\pi_1X$ (resp.~$\pi_1Y$).
Then $\rho(M,\al)=\rho(M',\al)$.

\begin{sloppypar}
In particular, if the image  of the restriction map
$\chi(\pi_1M,U(n))\to \chi(\pi_1X,U(n))$ lies in a path component of
$\chi(\pi_1X,U(n))$  and similarly for $Y$ (this holds e.g.\ when
$\chi(\pi_1X,U(n))$ and $\chi(\pi_1Y,U(n))$ are path connected),
then
$\rho(M,\al)=\rho(M',\al)$ for all $\al\in \chi(\pi_1M,U(n))$.
\end{sloppypar}
\end{thm}

\begin{proof} This follows by combining Theorems \ref {thm1.1} and \ref{FLWKL}.
\end{proof}

Notice that it is much more likely that the restrictions $\al|_X$ and $\al|_Y$
lie in the path component of the trivial connection than that $\al$
itself does,
since $\pi_1M$ is the free product of $\pi_1X$ and $\pi_1Y$ amalgamated over
$\pi_1\Si$. Hence in trying to deform $\al|_X$ and $\al|_Y$ in their
representation spaces one is no longer constrained by the relations imposed by
amalgamating over
$\pi_1\Si$.

As an application, if $X$ and $Y$ are manifolds with boundary which have  path
connected character varieties, and $f,g\co\del X \to \del Y$ are homotopic
diffeomorphisms, then $X\cup_f Y$ and $X\cup_g Y$ are homotopy equivalent and
$\rho(X\cup_f Y,\al)=\rho(X\cup_g Y,\al)$ for all $U(n)$ representations $\al$.
One can construct such examples so that $ X\cup_f Y$ and $X\cup_g Y$ are not
diffeomorphic.

The problem  of determining the number of path
components of $\chi(\pi,U(n))$ is  tricky.  Some examples of
groups with $\chi(\pi, U(n))$ path connected include $\pi$ free or  free
abelian. An interesting family of torsion-free groups with path connected
unitary representation  spaces are the 2-generator groups $\langle x,y, \ | \
x^p=y^q\rangle$ for $p,q$ relatively prime.  For a taste of the problem for
some 3-manifold groups
$\pi$ the reader might glance at
\cite{klassen-thesis} and \cite{KK-SF}.
Notice that   the isomorphism class of the
bundle $E\to X$ is fixed on any path component of $\chi(\pi_1X,U(n))$, and so
if $X$ admits non-isomorphic  flat bundles (e.g.\ if $H^2(X;\zz)$ contains
non-trivial torsion) then
$\chi(\pi_1X,U(n))$ cannot be path connected.

    Conjecture B  does not hold without the
   requirement that the homotopy equivalence behave  nicely on the
   boundary; this is exhibited by the example  at the end of Section
\ref{flwtheorem}.
  We will explore examples and applications of
Theorem
\ref{split} in a later  article.

\medskip
  We end this article with speculation concerning
the similarity between the constructions  of Section \ref{det bundles} and
the approach to studying TQFTs advocated in Atiyah's book \cite{atiyah-GPK}.

Take  $\Si$ to be a 2-manifold  and restrict  to
Riemannian metrics on $\Si$ which  have constant curvature $1,0,$ or $-1$.
Then $\mathscr{M}_\Si/\mathscr{D}^0_\Si= \mathscr{T}_\Si$ is the Teichm\"uller
space of $\Si$.     Thus  Theorem  \ref{thmcov} defines   a complex line
bundle with connection over
$\chi(\pi_1\Si,U(n))\times \mathscr{T}_\Si$, and given any 3-manifold $X$ with
boundary $\Si$ one obtains from $\rho(X,\al,g)$ a horizontal cross
section:
\begin{equation}\label{diag1}\begin{diagram}\dgmag{900}\dgsquash
\node[2]{\det(\ker A)^{-1}\otimes \det(\ker A_\theta) }\arrow{s}\\
\node{\chi(\pi_1X,U(n))\times
\mathscr{T}_\Si}\arrow{e,t}{i^*}\arrow{ne,t}{e^{\pi i
\rho}}
\node{\chi(\pi_1\Si,U(n))\times
\mathscr{T}_\Si}
\end{diagram}\end{equation}
In \cite{atiyah-GPK} a similar diagram is obtained: a determinant bundle over
$\chi(\pi_1\Si,U(n))$ is constructed as follows.  A metric $g\in
\mathscr{T}_\Si$  defines a holomorphic structure on $\Si$. This
defines a complex structure on $\chi(\pi_1\Si,U(n))$ by identifying it with the
moduli space of semi-stable holomorphic bundles over
$\Si$. Then one takes the
determinant bundle   $D\to \chi(\pi_1\Si,U(n))$ whose fiber over a point
corresponding to a holomorphic bundle is the determinant of the corresponding
$\bar{\del}$-operator.  The Quillen metric and the holomorphic structure on
this determinant bundle  determines a connection \cite{Quillen} which coincides
with $\nabla^Q$ \cite{Bismut-Freed}.  Viewing
holomorphic sections of this bundle as the fiber of a vector bundle over
$\mathscr{T}_\Si$ defines a bundle which was shown to admit a projectively
flat connection.

  Defining a horizontal cross-section of $D$
\[\begin{diagram}\dgmag{900}\dgsquash
\node[2]{D }\arrow{s}\\
\node{\chi(\pi_1X,U(n))}\arrow{e,t}{i^*}\arrow{ne}
\node{\chi(\pi_1\Si,U(n))}
\end{diagram}\]
for a 3-manifold $X$ with
boundary $\Si$ is  problematic from this point of view. An alternative set-up
is described in
  \cite{RSW}.  In that article a complex line bundle
$L\to \chi(\pi_1\Si,SU(2))$ is constructed using the Chern-Simons invariant
$cs$.  From the construction one  immediately  obtains the cross section
\[\begin{diagram}\dgmag{900}\dgsquash
\node[2]{L}\arrow{s}\\
\node{\chi(\pi_1X,SU(2))}\arrow{e,t}{i^*}\arrow{ne,t}{e^{2\pi i cs }}
\node{\chi(\pi_1\Si,SU(2))}
\end{diagram}\]
for any 3-manifold $X$ with boundary $\Si$ (this is a formalization of the fact
that on  a manifold with boundary, the Chern-Simons invariant is not a
$U(1)=\rr/\zz$ valued function but rather a cross section of a  $U(1)$
bundle over the moduli space of the boundary). This part of the construction is
independent of the choice of Riemannian metric on
$\Si$. However, endowing
$\Si$ with a metric defines a connection on $L$ for which the cross section is
horizontal.  It is shown in \cite{RSW} that the line bundles $L$ and $D$ are
isomorphic, linking the
two approaches.

On a closed 3-manifold, the $SU(n)$ Chern-Simons invariants and
$\rho$--invariants agree modulo $\zz$. This suggests that the set-up described in
\cite{atiyah-GPK,RSW} is related to our approach (encapsulated in the diagram
\eqref{diag1}). It would be an interesting project to establish a precise
relationship. Several complications arise. First, the $\rho$--invariant depends
on the metric on $\Si$, as we have established, but the Chern-Simons invariant
(viewed as a cross section of $L$) of a manifold with boundary is metric
independent. Second, it is not clear how to relate  the $\rho$ and
Chern-Simons invariants for a manifold with boundary. In the closed case the
Atiyah-Patodi-Singer index theorem provides the relationship, but generalizing
this argument would require an extension of the Atiyah-Patodi-Singer theorem
to manifolds with corners; a problem of significant interest that so far does
not have a complete solution.  Another interesting point is that the $\rho$--invariant
is  well--defined in $\rr$, not just $\rr/\zz$ as for the
Chern-Simons invariant, and the cut-and-paste formula \eqref{eq1.5} holds in
$\rr$.  This may give a clue as to how to refine the approach of
\cite{atiyah-GPK}.  Finally, since the set-up described in the present article
works in any odd dimension, it may provide a direction  to the problem of
constructing some TQFTs in higher dimensions.

\Addresses\recd
\end{document}